\documentclass[5p,number,sort&compress,times,fleqn, preprint]{elsarticle}

\usepackage{amsmath}
\usepackage{amssymb}
\usepackage{mathtools}

\usepackage[utf8]{inputenc}
\usepackage{subfig}
\usepackage{graphicx}

\usepackage[colorlinks, linktocpage=true, pdfpagemode=No,  linkcolor=blue, citecolor=blue, pagecolor=blue]{hyperref}

\usepackage{url}
\usepackage{tikz}
\usetikzlibrary{matrix}

\usepackage[noend]{algorithmic}
\usepackage{algorithm}

\newcounter{ALC@tempcntr}
\newcommand{\LCOMMENT}[1]{%
 	\setcounter{ALC@tempcntr}{\arabic{ALC@rem}}
     \setcounter{ALC@rem}{1}
     \item \textit{\small{// #1}}
     \setcounter{ALC@rem}{\arabic{ALC@tempcntr}}
 }%

\usepackage{csml}

\newcommand{\ADDEDTEXTMODET}[1]{\textcolor{csmlRed}{#1}}
\newcommand{\ADDEDT}[1]{
\ifmmode
  \text{\ADDEDTEXTMODET{$#1$}}
\else
  \ADDEDTEXTMODET{#1}
\fi}

\graphicspath{{./figs/}}

\begin{document}

\begin{frontmatter}

\title{Shape optimisation with multiresolution subdivision surfaces and immersed finite elements}

\author{Kosala Bandara}  
\author{Thomas R\"uberg}  
\author{Fehmi Cirak\corref{cor1}} 
\ead{f.cirak@eng.cam.ac.uk}

\cortext[cor1]{Corresponding author}

\address{Department of Engineering, University of Cambridge, Trumpington Street, Cambridge CB2 1PZ, U.K.}

\begin{abstract}
We develop a new optimisation technique that combines multiresolution subdivision surfaces for boundary description with immersed finite elements for the discretisation of the primal and adjoint problems of optimisation. Similar to wavelets multiresolution surfaces represent the domain boundary using a coarse control mesh and a sequence of detail vectors.  Based on the multiresolution decomposition efficient and fast algorithms are available for reconstructing control meshes of varying fineness.  During shape optimisation the vertex coordinates of control meshes are updated using the computed shape gradient information. By virtue of the multiresolution editing semantics, updating the coarse control mesh vertex coordinates leads to large-scale geometry changes and, conversely, updating the fine control mesh coordinates leads to small-scale geometry changes. In our computations we start by optimising the coarsest control mesh and refine it each time the cost function reaches a minimum. This approach effectively prevents the appearance of non-physical boundary geometry oscillations and control mesh pathologies, like inverted elements.
Independent of the fineness of the control mesh used for optimisation, on the immersed finite element grid the domain boundary is always represented with a relatively fine control mesh of fixed resolution. With the immersed finite element method there is no need to maintain an analysis suitable  domain mesh. In some of the presented two- and three-dimensional elasticity examples the topology derivative is used for creating new holes inside the domain. 
\end{abstract}


\end{frontmatter}

\newpage

\section{Introduction}

We consider the  shape optimisation of two- and three-dimen\-sional solids by combining multiresolution subdivision surfaces with immersed finite elements. As widely discussed in isogeometric analysis literature, the geometry representations used in today's computer aided design (CAD) and finite element analysis (FEA) software are inherently incompatible~\cite{Hughes:2005aa}. This is particularly limiting in shape optimisation during which a given CAD geometry model is to be iteratively updated based on the results of a finite element computation. The inherent shortcomings of present geometry and analysis representations have motivated the proliferation of various shape optimisation techniques. In the most prevalent approaches a surrogate geometry model~\cite{Braibant:1984aa, Haftka:1986aa,  Bletzinger:1991aa, olhoff1991cad, beux1994hierarchical, el2008multilevel, han2014adaptive} or the analysis mesh~\cite{Bletzinger:2010aa, le2011gradient} instead of the true CAD model is optimised, see also~\cite{bletzinger2014consistent} and references therein. Generally, it is tedious or impossible to map the optimised  surrogate geometry model or analysis mesh back to the original CAD model, which is essential for continuing with the design process and later for manufacturing purposes. Moreover, geometric design features are usually defined with respect to the CAD model and cannot be easily enforced on the surrogate model. Recently, the shape optimisation of shells, solids and other applications using isogeometric analysis has been explored; that is, through directly optimising the CAD geometry model~\cite{Cirak:2002aa,  kiendl2014isogeometric, Wall:2008aa, Bandara:2014aa}.

\begin{figure*}
	\centering 
	 \subfloat[][Original control mesh and  corresponding geometry.]{
         \includegraphics[width=0.33\linewidth]{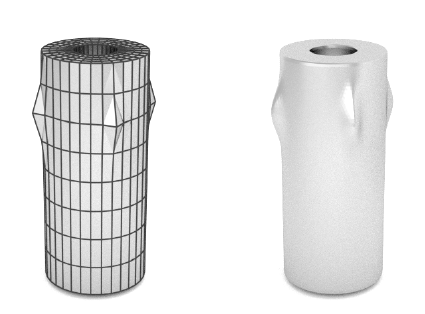}}
         \hspace{0.15\linewidth}
        \subfloat[][Edited control mesh and corresponding geometry.]{
            \includegraphics[width=0.33\linewidth]{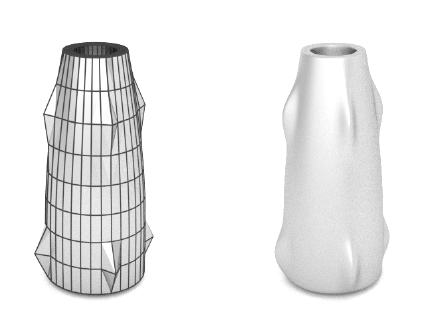}}
         \caption{Multiresolution editing of a cylindrical component with small scale geometric details in form of small bumps. The original geometry is shown on the left.  The edited version on the right contains both large and small scale changes.  The change of the overall shape from a cylinder to a cone constitutes a large scale change and  the additional bumps close to the bottom constitute small scale changes. The bumps  in the original geometry are automatically preserved during multiresolution editing.}
  \label{fig:editingIntroductory}
\end{figure*}

In the present work the domain boundary is represented with subdivision curves (in 2D) or  surfaces (in 3D). Although historically subdivision and related techniques have originated in computer graphics, they recently became available in several CAD software packages, including Autodesk Fusion 360, PTC Creo and CATIA. As will be demonstrated in this paper, subdivision curves/surfaces provide an elegant isogeometric, bidirectional mapping between the geometry and analysis models. In subdivision a geometry is described using a control mesh and a limiting process of repeated refinement~\cite{Zorin:2000aa, Peters:2008aa}.  The refinement rules are usually  adapted from knot refinement rules for b-splines~\cite{Lane:1980aa,Doo:1978aa, Catmull:1978aa}. The specific subdivision rules used in this work are derived from cubic b-splines. For surfaces we use the subdivision rules proposed by Catmull and Clark~\cite{ Catmull:1978aa}, which lead to smooth surfaces even in case of unstructured meshes with extraordinary vertices (i.e., domain vertices with number of adjacent edges different than four).  The hierarchy of control meshes underlying a subdivision surface lends itself naturally to multiresolution decomposition of geometry~\cite{zorin1997interactive, Lounsbery:1997aa}. To this end, suitable operators are available for decomposing a geometry in  a coarse control mesh and a sequence of detail vectors similar to wavelets. Subsequently, it becomes possible to reconstruct on-the-fly control meshes of any fineness and to edit their vertex positions.  The size of the geometric region influenced by each vertex depends on the resolution of the control mesh, editing coarser levels leads to large-scale changes while editing finer levels lead to small-scale changes. Importantly, after applying large-scale changes to the limit geometry the available  finer detail vectors can be automatically added to obtain the new limit geometry, cf. Figure~\ref{fig:editingIntroductory}.

During gradient-based shape optimisation the \emph{shape gradient} is computed to determine the boundary perturbations which lead to a reduction in a given cost function. We consider the adjoint approach as applied to the continuous shape optimisation problem for computing the shape gradient~\cite{Haug:1986aa, Sokolowski:1992aa, Haslinger:2003aa}. The shape gradient is a function of the solution of the original (primal) and a complimentary adjoint boundary value problem.   In case of compliance minimisation the primal and adjoint solutions are up to the sign identical and, hence, only the primal boundary value problem has to be solved.  We discretise the primal and adjoint  problems with immersed, or embedded,  finite elements, see, e.g.,~\cite{Ruberg:2010aa, Sanches:2011aa, Schillinger:2012aa, Ruberg:2011aa}, which have clear advantages when applied to structural shape optimisation \cite{Haslinger:2003aa, Norato:2004aa, Allaire:2004aa}. The geometry of the domain boundary can be updated without needing to generate a new, or to smooth, an existing domain mesh.  On the fixed non-boundary-conforming immersed finite  element grid, we use tensor-product b-splines as basis functions and enforce Dirichlet boundary conditions with the Nitsche technique. 

As mentioned, the developed multiresolution optimisation approach relies on multiresolution subdivision curves/surfaces for geometry description and the immersed finite element method for domain discretisation. 
 The multiresolution representation of the domain boundaries allows us to describe the same geometry with  control meshes of different resolution for analysis and optimisation purposes.  For finite element analysis a relatively fine control mesh is used in order to faithfully describe the domain boundaries on the immersed  grid. In contrast, the degrees of freedom in optimisation (i.e., design variables) are chosen as the vertex coordinates of a coarser control mesh. Usually, better shapes can be found by starting with a coarse control mesh and optimising increasingly finer control meshes. During the optimisation iterations the refinement level of the control mesh is increased each time a minimum is reached. 
Informally, the superior performance of the multiresolution approach can be explained with the correlation between the number of design variables and the number of local minima of the cost function. Having initially fewer local minima reduces the possibility of landing in a local minimum. 
It bears emphasis that with the employed,  wavelet-like multiresolution decomposition, the coarse control mesh for optimisation can be chosen independently of the size of the features present in the geometry. For instance, the control mesh for optimisation can be chosen much coarser than any fillets or small-scale surface undulations present on the to be optimised geometry. 

 There are a number of prior works, especially in aerodynamic shape optimisation, that use hierarchical and adaptive geometry representations, see e.g.~\cite{Borzi2013} and references therein.  For instance, in~\cite{el2008multilevel} Bezier basis functions and degree elevation and in~\cite{han2014adaptive}  b-splines and knot insertion  are considered to create a hierarchy of geometry representations. Most of these papers primary aim to speed up the optimisation process by reducing the number of optimisation variables or by employing multigrid techniques.  However, the added benefit  of hierarchical representations in reducing  the parameterisation dependency of the  final results is also often noted. In comparison to the mentioned techniques, the advantages of multiresolution subdivision surfaces are: (i) the ability to represent geometries with arbitrary topology, (ii) wavelet-like multiresolution representation of the geometry,  and (iii) the ease of integration with CAD packages that use subdivision.

This paper is organised as follows. Section~\ref{sec:governing} introduces the governing equations for linear elastic shape optimisation.  Specifically, the  derivation of the continuous shape gradient using the adjoint approach is illustrated. Subsequently, Section~\ref{sec:immersedFE} discusses the discretisation of the primary and the adjoint boundary value problems with an immersed finite element method. Section~\ref{sec:multires} forms the core of the paper and introduces multiresolution optimisation. First the derivation of the cubic b-spline and Catmull-Clark subdivision rules from the well known b-spline refinement relations is demonstrated. After that, the multiresolution decomposition of subdivision surfaces and its use for multiresolution editing are shown. Upon this basis we introduce our multiresolution shape optimisation algorithm. Finally, Section~\ref{sec:examples}  presents several two- and three-dimensional examples of increasing complexity to demonstrate the efficiency and robustness of multiresolution optimisation. In some of the examples the topology derivative is considered to introduce new holes in the domain.

\section{Governing equations \label{sec:governing}}
%
We  introduce in this section the shape optimisation of two- and three-dimen\-sional linear elastic solids. In addition to the linear  elasticity boundary value problem a cost function that is to be minimised is considered.  For computing the shape derivatives, that is the derivatives of the cost function with respect to the domain perturbations, we chose to use the analytic adjoint approach. As it will become clear in Section~\ref{sec:immersedFE}, with the analytic approach it is straightforward to compute the shape gradients using the immersed finite element technique.
In the following, we  briefly review few key results from shape calculus. See, for instance, the monographs~\cite{Haug:1986aa, Sokolowski:1992aa,Delfour:2001aa,Haslinger:2003aa} for a more detailed discussion.

\subsection{Linear elasticity \label{sec:linElasticity}}
%
The equilibrium equation for a solid body with the domain~$\Omega$ is given by
\begin{subequations} \label{eq:elasticityStrong}
  \begin{align}
  		\nabla \cdot \vec \sigma (\vec u)+ \vec f &= \vec 0 &  \text {in } & \Omega \, ,  \\
		\vec u &= \vec 0 & \text {on } & \Gamma_D \, , \\
		\vec \sigma (\vec u) \vec n &= \overline{\vec t} & \text{on } & \Gamma_N \,  ,
 	\end{align}	
\end{subequations}
where $\vec \sigma$ is the stress tensor, $\vec u$ is the displacement vector,  $\vec f$ is the external load vector and $ \overline{ \vec t}$ is the prescribed traction on the Neumann boundary $\Gamma_N$ with the outward normal $\vec n$. On the Dirichlet boundary $\Gamma_D$, for simplicity,  only homogenous  boundary conditions are prescribed.

We assume a homogenous linear elastic material model 
\begin{equation}
	\vec \sigma (\vec u) = \vec C: \vec \epsilon (\vec u)	
\end{equation}
with the fourth order constitutive tensor $\vec C$  and linear elastic strain tensor 
\begin{equation}
	\vec \epsilon (\vec u) = \frac{1}{2} (\nabla \vec u + \nabla^\trans \vec u)  \, .
\end{equation}
%

\subsection{Shape optimisation}
%

\subsubsection{Formulation}
%

The cost function to be minimised 
\begin{equation} \label{eq:costFunction}
	 J (\Omega, \vec u) \rightarrow \min
\end{equation}
depends on the to be optimised domain $\Omega$ and the unknown solution~$\vec u$. The most common cost functions are integrals over the domain~$\Omega$ or its boundary~$\Gamma$.  For brevity and without loss of generality, we consider in the following only the structural compliance as the cost function 
\begin{align}\label{eq:compliance}
	\begin{split}
	  J (\Omega, \vec u )  &= \int_{\Omega} \vec \sigma (\vec u ) : \vec \epsilon (\vec u) \,  \dif \Omega   \\ 
						& = \int_{\Omega} \vec f \cdot \vec u \,  \dif \Omega  +  \int_{\Gamma_N}  \overline{\vec t} \cdot \vec u \,  \dif \Gamma \, .
	\end{split}
\end{align}	
It is straightforward to adapt the subsequent derivations to other common cost functions, such as integrals of stresses or displacements.  

The minimisation of the cost function  $J (\Omega, u)$ with the boundary value problem~\eqref{eq:elasticityStrong}  as a constraint can be achieved by means of the Lagrangian function
\begin{align} \label{eq:lagrangianWeak}
	\begin{split}
		  L (\Omega, \vec u, \vec \lambda)  &=   J (\Omega, \vec u ) 
		\\ &+ \int _{\Omega}    \nabla \vec \lambda :  \vec \sigma (\vec u) \, \dif \Omega   -  \int _{\Omega} \vec \lambda \cdot \vec f   \, \dif \Omega    
		\\& - \int_{\Gamma_D}  \vec u \cdot (  \vec C : \nabla \vec \lambda )   \vec n + \vec \lambda \cdot \vec \sigma (\vec u)  \vec n \, \dif \Gamma  
		\\& -  \int_{\Gamma_N } \vec \lambda  \cdot   \overline {\vec t}  \, \dif \Gamma \, ,
	\end{split}	
\end{align}
where $\vec \lambda$ is a Lagrange parameter. $  L (\Omega, \vec u, \vec \lambda) $ depends on the unknown domain $\Omega$, the displacement field $\vec u$ and the Lagrange parameter $\vec \lambda$. 
The stationarity condition for the Lagrangian~\eqref{eq:lagrangianWeak}, i.e. $\delta   L (\Omega, \vec u, \vec \lambda)=0 $, 
provides the complete set of equations that describe the shape optimisation problem.  In the following we consider one after the other the variation of   $   L (\Omega, \vec u, \vec \lambda) $ with respect to the Lagrange parameter~$\vec \lambda$,  displacements~$\vec u$ and domain~$\Omega$.

\subsubsection{Primal problem}
%
 The  variation of  the Lagrangian $  L (\Omega, \vec u, \vec \lambda) $ with respect to $\vec \lambda$  reads
\begin{align} \label{eq:varWrtLambda}
	 \frac{\partial   L }{\partial  \vec \lambda} \delta  \vec \lambda   +  \frac{\partial   L }{\partial ( \nabla \vec \lambda)} \delta ( \nabla \vec \lambda ) \notag = & - \int _{\Omega} \delta \vec \lambda \cdot  [  \nabla \cdot \vec \sigma (\vec u)+ \vec f  ] \, \dif \Omega 
  \\ & - \int_{\Gamma_D}  \vec u \cdot  \left ( \vec C :   \nabla ( \delta  \vec \lambda)  \right ) \vec n \, \dif \Gamma 
  \\ & - \int_{\Gamma_N}   \delta \vec \lambda \cdot [ \overline {\vec t} - \vec \sigma (\vec u) \vec n ] \, \dif \Gamma = 0 \notag \,  , 
\end{align} 
where we used the divergence theorem and $\delta (\nabla \vec \lambda) = \nabla (\delta \vec \lambda)$.   For arbitrary $\delta \vec \lambda$  it is evident that  \eqref{eq:varWrtLambda} is equivalent to the linear elasticity equations~\eqref{eq:elasticityStrong}. 

\subsubsection{Adjoint problem}
%
Next, we consider the variation of the Lagrangian $  L (\Omega, \vec u, \vec \lambda) $  with respect to the displacements $\vec u$
\begin{align}
	 \frac{\partial   L }{\partial  \vec u} \delta  \vec u    &+  \frac{\partial   L }{\partial ( \nabla \vec u)} \delta ( \nabla \vec u ) \notag =  
	 \frac{\partial J}{\partial \vec u} \delta \vec u
	\\ &   +\int_{\Omega} \nabla \vec \lambda :   \vec C : \delta (\nabla \vec u)  \, \dif  \Omega 
	 \\& - \int_{\Gamma_D} \delta \vec u \cdot (\vec C : \nabla \vec \lambda ) \vec n + \vec \lambda \cdot \left (  \vec C :  \delta(\nabla \vec u) \right )  \vec n \, \dif \Gamma= 0 \notag \, .
\end{align}
After introducing the cost function~\eqref{eq:compliance}  and reformulating the domain term with the divergence theorem we obtain 
\begin{align} 
	\begin{split}
	 & \int_\Omega   \vec f   \cdot \delta \vec u  \, \dif \Omega   + \int_{\Gamma_N} \overline{\vec t} \cdot \delta \vec u \, \dif \Gamma
	  	  -\int_{\Omega} \delta \vec u \cdot (  \nabla \cdot \vec \sigma (\vec \lambda) ) \, \dif  \Omega  
	 \\&  - \int_{\Gamma_D} \vec \lambda \cdot \left (  \vec C :   \nabla ( \delta \vec u) \right )  \vec n   \, \dif \Gamma + \int_{\Gamma_N}  \delta \vec u \cdot  \vec \sigma (\vec \lambda)  \vec n  \, \dif \Gamma= 0  \, . 
	 \end{split}
\end{align}
The corresponding boundary value problem, referred to as the adjoint problem, reads 
\begin{subequations} \label{eq:adjointStrong}
  \begin{align}
  		 \nabla \cdot \vec \sigma (\vec \lambda) - \vec f  &= \vec 0 &  \text {in } & \Omega \, , \\
		\vec \lambda &= \vec 0 & \text {on } & \Gamma_D \, , \\
		\vec \sigma (\vec \lambda) \vec n &=  - \overline{\vec t} & \text{on } & \Gamma_N \, .
 	\end{align}	
\end{subequations}
By comparing this adjoint problem with the primal problem~\eqref{eq:elasticityStrong}  we deduce that $\vec \lambda = - \vec u$. Note that this holds only when the cost function is the structural compliance~\eqref{eq:compliance}.  

\subsubsection{Shape derivative}
%

\begin{figure*}
	\centering
	\includegraphics[scale=0.9]{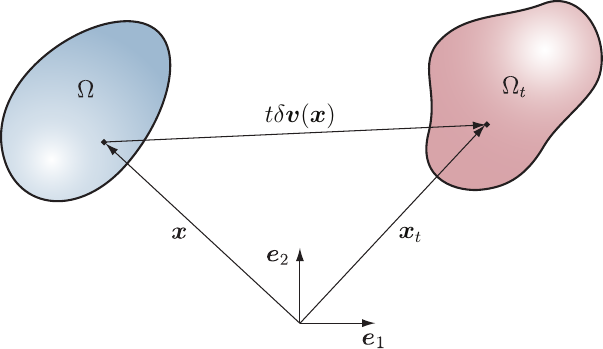}
	\caption{Reference  and the perturbed domains (left and right, respectively). \label{fig:domains} }
\end{figure*}

Finally, we consider the  variation of $  L (\Omega, \vec u, \vec \lambda)$  with respect to the problem domain~$\Omega$, which is also referred to as the shape derivative. To this end, we first define  a linear mapping which   maps a given domain  $\Omega$ into a perturbed domain $\Omega_t$, see Figure~\ref{fig:domains}. With this mapping  a material point with the coordinate $\vec x \in \Omega$ is mapped onto  
\begin{equation} \label{eq:perturbedDomain}
	\vec x_t = \vec {x} +t  \delta \vec v  \,  , 
\end{equation}
where    $\delta \vec v$ is a prescribed vector field and $t$ is a scalar parameter.  
In the usual  continuum mechanics terminology, $\Omega$ is the reference configuration, $\Omega_t$ is the current configuration, $\delta \vec v$ is the material velocity vector and $t$ is the (pseudo-) time. In shape optimisation literature the mapping~\eqref{eq:perturbedDomain} is usually expressed as
\begin{equation}
	 \Omega_t =  \Omega + t \delta \vec v  \, .
\end{equation}

Before attempting  the variation of  $  L (\Omega, \vec u, \vec \lambda) $, we first give the variation of generic volume and surfaces integrals with respect  to $\delta \vec v$.  The variation of a domain  integral  of  the scalar function  $\psi(\vec x_t)$ 
\begin{equation}
	  I_1  (\Omega) = \int_{\Omega}   \psi (\vec x ) \, \dif \Omega
\end{equation}
at the reference configuration $\vec \Omega$ in the direction of $\delta \vec v$ is defined as 
\begin{equation} \label{eq:volDerivShape}
	\dfrac{\partial   I_1}{\partial \Omega}  \delta \vec v =   \frac{d}{dt}   I_1(\vec x + t \delta \vec v ) \Big |_{t=0} = \frac{d}{dt} \int_{\Omega_t}  \psi ({\vec x_t}) \, \dif \Omega_t \Big |_{t=0}  \, .
\end{equation}
Transforming this integral into the reference configuration yields  
\begin{equation} \label{eq:reynoldsLike} 
\begin{split}
	\dfrac{\partial   I_1}{\partial \Omega} \delta \vec v &=  
	 \frac{d}{dt} \int_{\Omega}    \psi (\vec x_t)  j (\vec x_t)  \, \dif \Omega \Big |_{t=0} 
	 \\ &=   \int_{\Omega}   \left ( \frac{d}{dt} \psi (\vec x_t) j (\vec x_t)+ 	\psi (\vec x_t) \frac{d}{dt} j (\vec x_t)  \right ) \, \dif \Omega \Big |_{t=0}  \, 
	\end{split}
\end{equation}
with the determinant of the mapping
\begin{equation*}
	j (\vec x_t) = \det \frac{\partial \vec x_t}{\partial \vec x}  \, ,
\end{equation*}
which has, according to, e.g.,~\cite{Holzapfel:2000aa, Delfour:2001aa},  the derivative
\begin{equation} \label{eq:detDer}
	 \frac{d}{dt}  j(\vec x_t) \Big |_ {t=0}=
	  \nabla \cdot \delta \vec v \, .
\end{equation}
After introducing~\eqref{eq:detDer} into~\eqref{eq:reynoldsLike} and applying the divergence theorem we obtain
\begin{equation} \label{eq:reynolds} 
\begin{split}
	\dfrac{\partial I_1}{\partial \Omega} \delta \vec v &=  
	   \int_{\Omega}   \left ( \nabla \psi (\vec x) \cdot   \delta \vec v + 	\psi (\vec x) \nabla \cdot \delta \vec v  \right ) \, \dif \Omega
	 \\ & =   \int_{\Gamma} \psi (\vec x) (\delta \vec v \cdot \vec n) \, \dif \Gamma \, , 
	\end{split}
\end{equation}
where $\vec n$ is the unit normal to the boundary. Notice that this integral is zero when the perturbation direction $\delta \vec v$ is chosen tangential to the boundary.  Perturbations tangential to the boundary do not lead to a change in $I_1$.

Although not used in this work,  we give for completeness  the variation of the boundary integral of the scalar function~$\psi (\vec  x_t)$,  
\begin{equation}
		  I_2  (\Gamma) = \int_{\Gamma}   \psi (\vec x) \, \dif \Gamma \, .
\end{equation}
The variation of this integral at the reference configuration $\vec x$ in the direction of $\delta \vec v$ reads
\begin{align} \label{eq:surfaceReynolds}
		\dfrac{ \partial   I_2  }{\partial \Gamma} \delta \vec v =   \int_{\Gamma} \left  ( \nabla \psi  ( \vec x) \cdot \vec n  + H (\vec x) \psi (\vec x) \right ) (\delta \vec v \cdot \vec n) \, \dif \Gamma \, ,
\end{align}
where $H(\vec x)$ is the mean curvature on $\Gamma$, see~\cite{Sokolowski:1992aa, Delfour:2001aa}. 

We can now write the variation of the Lagrangian  $  L (\Omega, \vec u, \vec \lambda) $ in the direction of $\delta \vec v$ using the results~\eqref{eq:reynolds} and~\eqref{eq:surfaceReynolds}.  For practical shape optimisation problems in solid mechanics we usually have only boundary variations of the form 
\begin{align} \label{eq:variationHat}
	\begin{split}
	 	\delta {\hat{\vec v} } &= \vec 0 \quad \text{on } \;  \Gamma_D \, ,\\
		\delta  {\hat{\vec v}} &= \vec 0 \quad \text {on } \;  \Gamma_N  \; \text  { with }  \; \vec \sigma \vec n = \overline{ \vec t  }\, , \\
		\delta  {\hat{\vec v}} &\neq \vec 0 \quad \text {on } \;  \Gamma_N  \; \text  { with }  \; \vec \sigma \vec n =  \vec 0  \, .
	\end{split}
\end{align}
This means that only parts of the Neumann boundary $\Gamma_N$  with no traction are free to move during shape optimisation. With the result~\eqref{eq:reynolds} at hand the variation of the Lagrangian~\eqref{eq:lagrangianWeak} in the direction~$\delta \vec{\hat v}$  reads
\begin{equation}
	\begin{split}
	\frac{\partial   L}{\partial \Omega} \delta \hat {\vec v} &= \frac{\partial   J(\Omega, \vec u) }{\partial \Omega} \delta \hat { \vec  v}  \\ & + 
		\int_{\Gamma_N} \left ( \nabla \vec \lambda : \vec \sigma (\vec u) - \vec \lambda \cdot \vec f ) \right (\delta \hat{\vec v} \cdot \vec n )\, \dif \Gamma \, .
	\end{split}
\end{equation}
For structural compliance~\eqref{eq:compliance} as the cost function (i.e., $\vec \lambda = - \vec u$) we obtain 
\begin{align}~\label{eq:shapeDerivFinal}
	\begin{split}
		\frac{\partial   L}{\partial \Omega} \delta \hat {\vec v} & =
		\int_{\Gamma_N} \left ( 2 \vec u \cdot \vec f   -  \nabla \vec u : \vec \sigma (\vec u)   \right ) (\delta \hat{\vec v} \cdot \vec n )\, \dif \Gamma  \\
		& = \int_{\Gamma_N} g (\vec u) (\delta \hat{\vec v} \cdot \vec n )\, \dif \Gamma  \, , 
	\end{split}	
\end{align}
where $g(\vec u)$  is  the shape kernel function. It is worth emphasising that without restricting $\delta \hat { \vec v}$ as stated in~\eqref{eq:variationHat}, the shape derivate would contain several more terms. Moreover, for cost functions other than the structural compliance, the kernel function usually is also dependent on the adjoint solution $\vec \lambda$.

During the iterative shape optimisation the shape kernel function $g(\vec u)$ is used as gradient information. In order to achieve a maximum decrease in the cost function the boundary perturbation is chosen proportional to 
\begin{equation} ~\label{eq:shapeDerivKernel}
	\delta \hat{\vec v}     = -  g ( \vec u )  \vec n  
\end{equation}
such that 
\begin{equation}
	\frac{\partial   L}{\partial \Omega} \delta \hat {\vec v} =
		- \int_{\Gamma_N}      g ( \vec u )^2   \, \dif \Gamma \, .
\end{equation}
%

\section{Immersed finite element discretization \label{sec:immersedFE}}
%
The shape derivatives introduced in the previous section depend on the solution of the primal and adjoint problems~\eqref{eq:elasticityStrong} and~\eqref{eq:adjointStrong}, respectively. Although for compliance optimisation the primal and adjoint solutions are identical (up to sign),  this is not generally the case for other cost functions.  During the iterative shape optimisation, both boundary value problems have to be repeatedly solved on constantly evolving domains.  In a conventional finite element setting this requires frequent mesh smoothing or  updating. Therefore, immersed,  or embedded, grid finite element approaches that do not require remeshing  have clear advantages in shape optimisation~\cite{Haslinger:2003aa, Norato:2004aa, van2013level}. In the present work, we use an immersed finite element technique that we previously developed in the context of incompressible fluid-structure interaction~\cite{Ruberg:2011aa,Ruberg:2014aa}. The key features of which are: (i) the weak enforcement of Dirichlet boundary conditions with the Nitsche technique, (ii) the use of isoparametric b-spline basis  functions for discretisation, and (iii) the numerically robust boundary and cut-cell treatment. In the following we provide a brief summary of our discretisation method. Although we only discuss the discretization of the primal problem~\eqref{eq:elasticityStrong}, the same derivations also apply to the adjoint problem~\eqref{eq:adjointStrong}.

\subsection{Weak form of the equilibrium equations}
%
%
For the  linear elastic solid introduced in Section~\ref{sec:linElasticity}, the weak form of the equilibrium equations~\eqref{eq:elasticityStrong} reads
\begin{align}
  \label{eq:weakForm}
  \begin{split}
	  \int_\Omega \vec{\sigma}(\vec u) : \vec{\epsilon}(\delta \vec u) \D{\Omega} & =
	  \int_\Omega \vec{f}\cdot \delta \vec u \D{\Omega}  + \int_{\Gamma_N} \overline{\vec t} \cdot
	  \delta \vec u \D{\Gamma} \\ & + \int_{\Gamma_D} \vec t(\vec u) \cdot \delta \vec u \D{\Gamma}   \, ,
	\end{split}  
\end{align}
where  $\delta \vec u$ are test functions, which are here assumed not to be zero on the Dirichlet boundary. This assumption is necessary because we use non-boundary-fitting meshes. The weak form~\eqref{eq:weakForm} is not coercive and would lead, for instance, to a singular system of equations when discretised. In order to render~\eqref{eq:weakForm} coercive we use the  consistent penalty method proposed by Nitsche~\cite{Nitsche:1971aa}, which leads to
\begin{align}
  \label{eq:weakFormNitsche}
  \begin{split}
	  & \int_\Omega \vec{\sigma}(\vec u) : \vec{\epsilon}(\delta \vec u) \D{\Omega}   
	   - \int_{\Gamma_D}\vec t(\vec u) \cdot \delta \vec u \D{\Gamma} \\  - &
  	\int_{\Gamma_D} \vec u \cdot \vec t(\delta \vec u) \D{\Gamma}  +  \frac{\gamma}{h} \int_{\Gamma_D}
  	\vec u \cdot \delta \vec u \D{\Gamma}  
  	\\  = & \int_\Omega \vec{f}\cdot \delta \vec u \D{\Omega} + \int_{\Gamma_N} \overline{\vec t} \cdot  \delta \vec u \D{\Gamma}  
  \end{split}
\end{align}
with the penalty parameter~$\gamma > 0$ and the characteristic finite element size~$h$. In contrast to conventional penalty methods, the parameter $\gamma$ in the Nitsche method is only required for numerical stability and typically a small value is sufficient.

\subsection{Finite element discretisation}
\begin{figure}
	\centering 
        \includegraphics[width=0.7\linewidth]{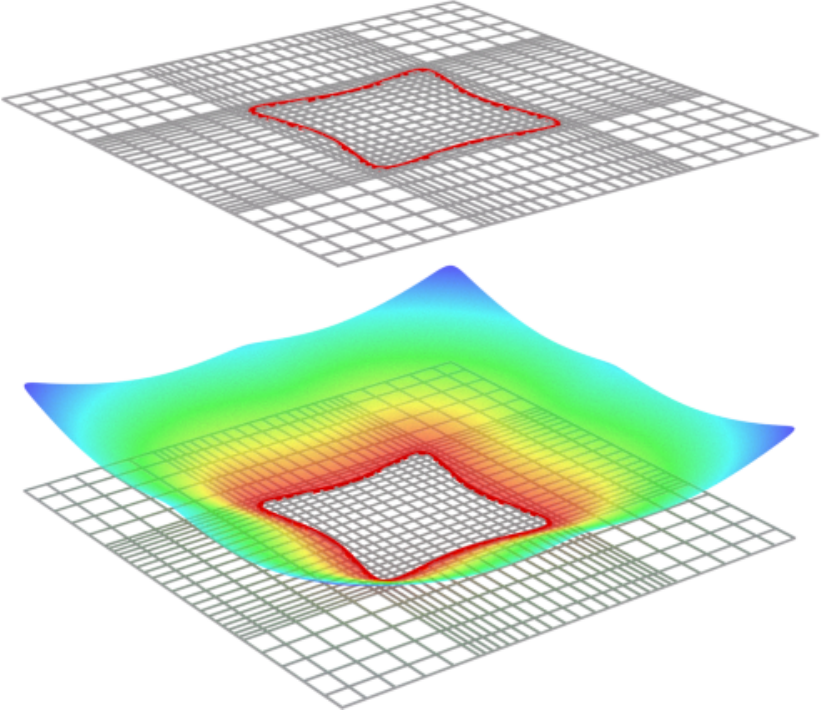}	
 	\caption{Domain discretization with the immersed finite element method method. The top figures shows a (logically) Cartesian grid and a red spline curve describing the domain boundary. On the Cartesian grid the boundary is represented using the signed distance function shown below. The zero isocontour of the distance function provides an approximation to the original spline curve.}
  \label{fig:immersedDomain}
\end{figure}

%
%
We use a logically Cartesian grid and the associated tensor-product b-spline basis functions for discretizing the weak form~\eqref{eq:weakFormNitsche}.  The grid has to have the connectivity of a Cartesian grid but the cell sizes need not to be uniform. Figure~\ref{fig:immersedDomain} shows a typical setup in two space dimensions.  The Cartesian grid facilitates the use of tensor-product b-spline basis functions, which have a number of appealing properties known from isogeometric analysis. Specifically, in shape optimisation the smoothness of higher-order b-splines leads to a shape kernel function~\eqref{eq:shapeDerivFinal}, which is continuous across element boundaries. This leads  to optimisation algorithms that are more robust than the ones based on $C^0$-continuous shape functions and discontinuous shape gradients. 

According to the isoparametric concept, we approximate the domain geometry and the solution with tensor-product b-splines.     
\begin{equation}
  \label{eq:feApprox}
   \vec{x} \approx  \vec{x}_h(\vec{\xi}) = \sum_{\vec{i}}  B^\alpha_{\vec{i}}(\vec{\xi}) \vec{x}_{\vec{i}} \, , \quad 		\vec{u} \approx \vec{u}_{h}(\vec{\xi}) =  \sum_{\vec{i}} B_{\vec{i}}^{\alpha}(\vec{\xi})   \vec{u}_{\vec{i}} \, ,
\end{equation}
where~$\vec{\xi}$ are the parametric coordinates in the knot space and~$\vec i$ is the control point index with the corresponding b-spline~$B^\alpha_{\vec{ i}}$  of degree $\alpha$. In optimisation computations,  we usually use quadratic b-splines with $\alpha=2$, which are $C^1$-continuous and lead to continuous shape gradients. The control point coefficients~$\vec{x}_{\vec{i}}$ and~$ \vec{u}_{\vec{i}}$ are the nodal coordinates and displacements, respectively. The discrete finite element equations are obtained  by introducing the approximation equations~\eqref{eq:feApprox} into the weak form~\eqref{eq:weakFormNitsche} and subsequent element-by-element numerical integration. The elements which are only partially covered by the solid, so-called cut-cells,  are first triangulated prior to integration, see~\cite{Ruberg:2014aa}. The ill-conditioning of the system matrix due to small cut-cells is avoided with an extension approach originally proposed in~\cite{hoellig}. Finally, note that the shape gradient in the cut-cells is computed by simply evaluating~\eqref{eq:shapeDerivFinal}. The advantage of the analytic adjoint formulation here is that the derivatives of the cut-cells with respect to the boundary position are not needed.  

On the logically Cartesian grid we represent the domain boundaries implicitly with a signed distance function (or, in other terms, a level set). This is done despite the fact that we have a parametric representation of the  domain in form of a multiresolution subdivision surface, see Section~\ref{sec:multires}. The aim of the switch from a parametric to an implicit representation is to eliminate pathological geometries and topologies, like multiple crossings of a cell edge by the boundary. For  optimisation problems with large boundary deformations and topology changes the  low-pass  filtering of the geometry provided through the switch to an implicit representation makes the  finite element computations more robust.

\section{Multiresolution optimisation \label{sec:multires}}
%
On the logically Cartesian discretisation grid we represent the domain boundaries  either with subdivision curves or subdivision surfaces, depending on the dimensionality of the problem. The inherent hierarchy of subdivision schemes lends itself to  multiresolution representation and editing of curves and surfaces. The specific subdivision scheme that we use yields in  the curve case cubic b-splines and in the surface case cubic tensor-product b-splines~\cite{Catmull:1978aa}. Subdivision surfaces yield smooth surfaces even for unstructured surface meshes with non-tensor product structure. We refer to the monograph~\cite{Zorin:2000aa} as an introduction to subdivision surfaces and multiresolution editing in geometric modelling and animation. 
%

\subsection{Subdivision scheme for one-dimensional cubic b-splines}
%
The refinability property of cubic b-splines can be utilised to derive a corresponding subdivision scheme.  To illustrate this, we consider the coarse knot sequence $\xi_i= 0, \,1,  \,2,  \, 3, \ldots$ and the corresponding fine knot sequence $\tilde \xi_i= 0, \,  0.5, \, 1, \, 1.5,   \, 2,  \, 2.5, \, 3,   \ldots$.  We denote the b-splines on the coarse knot sequence with $B_i (\xi)$  and the ones on the fine knot sequence with $\tilde B_i ( \xi)$.  According to the b-spline refinability equation, see, e.g.,~\cite{Zorin:2000aa, Prautzsch:2002aa},  it is possible to represent the  coarse b-splines as a linear combination of the fine b-splines
\begin{equation} \label{eq:coarseFineBsplines}
	B_i (\xi) = \sum_j S_{ij} \tilde B_j (\xi) \, ,
\end{equation}
where $S_{ij}$ is the subdivision matrix with the components

\begin{figure}
	\centering 
        \includegraphics[scale=0.9]{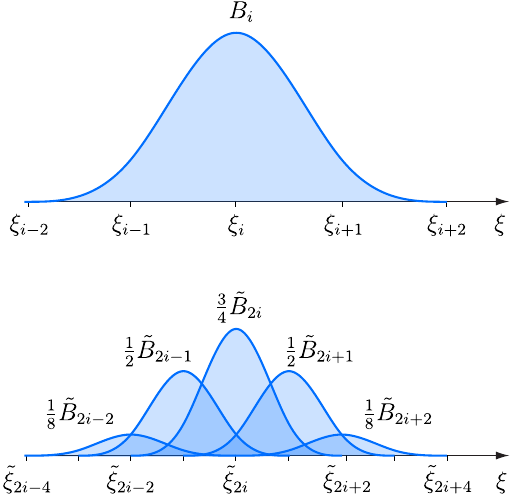}	
 	\caption{Refinement relation for cubic b-splines.}
  \label{fig:refinementCubicSpline}
\end{figure}

\setcounter{MaxMatrixCols}{12}
\begin{equation} \label{eq:subdivMatrix}
S_{ij} = 
\begin{bmatrix}
\cdots 	& \cdots 			& \cdots 			&  \cdots 			& \cdots 			& \cdots  			& \cdots  			& \cdots      		 & \cdots	  \\[0.5em]
\cdots 	& \dfrac{1}{2} 	& \dfrac{1}{8}	&  0 					& 0 					& 0  					& 0  					& 0  					 & \cdots	 \\[1em]
\cdots   & \dfrac{1}{2} 	& \dfrac{3}{4} 	& \dfrac{1}{2} 	& \dfrac{1}{8} 	& 0 					& 0 					& 0		 		      & \cdots   	 \\[1em]
\cdots   & 0 					& \dfrac{1}{8} 	& \dfrac{1}{2} 	& \dfrac{3}{4} 	& \dfrac{1}{2} 	& \dfrac{1}{8} 	& 0 		     	  	 & \cdots 	 \\[1em]
\cdots   & 0 					& 0 					& 0 					& \dfrac{1}{8} 	& \dfrac{1}{2} 	& \dfrac{3}{4} 	& \dfrac{1}{2}     & \cdots    \\[1em]
\cdots   & 0 					& 0 					& 0 					& 0 					& 0 					& \dfrac{1}{8} 	& \dfrac{1}{2} 	 & \cdots    \\[0.5em]
\cdots 	& \cdots 			& \cdots 			& \cdots 			& \cdots 			& \cdots 			& \cdots 			 & \cdots 			 & \cdots 
\end{bmatrix} \, .
\end{equation}
Each row of this banded matrix has the same five non-zero components, always shifted by two columns relative two adjacent rows. As shown in Figure~\ref{fig:refinementCubicSpline} each row expresses how a coarse b-spline can be obtained as the weighted sum of fine b-splines. It is evident that the exact structure and components of the matrix~\eqref{eq:subdivMatrix} depends on the degree of the considered b-splines. The components of the subdivision matrix for different polynomial degrees can be found, e.g., in~\cite{Ruberg:2011aa}.

\begin{figure}
	\centering 
	\subfloat[][Even vertex] {
  		\includegraphics[scale=1]{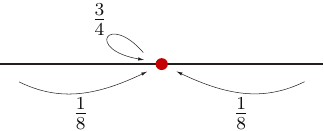}} \\ 
  	\subfloat[][Odd vertex]{
		\includegraphics[scale=1]{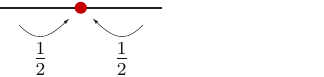}}
	  	\caption{Subdivision stencils for cubic b-splines.}
  \label{fig:cubicMasks}
\end{figure}

Next, we consider a spline curve defined in terms of the coarse b-splines  and the corresponding control vertices with the coordinates $ \vec{x}_{i} $, i.e., 
\begin{equation} \label{eq:coarseSpline}
	  \vec{x}_h(\xi) = \sum_{i}  B_i(\xi) \vec{x}_i \, .
\end{equation}
The control polygon of the spline is obtained by linearly connecting the control points $\vec{x}_i$.   Because of the refinability property \eqref{eq:coarseFineBsplines}, the spline curve~\eqref{eq:coarseSpline} can also be represented with the fine b-splines $\tilde B_j (\xi)$. To show this, we  introduce the refinement relation~\eqref{eq:coarseFineBsplines} into~\eqref{eq:coarseSpline} 
\begin{equation}
	 \vec{x}_h(\xi) = \sum_{i}    \sum_j S_{ij} \tilde B_j (\xi)  \vec{x}_{i}  =   \sum_j  \tilde B_j (\xi)   \left (  \sum_{i}  S_{ij}  \vec{x}_{i}  \right ) \, .
\end{equation}
Hence, choosing the fine control vertex coordinates with
\begin{equation} \label{eq:cpRefinement}
	 \tilde{\vec x}_j= \sum_i S_{ij}  \vec{x}_{i}
\end{equation}
ensures that both the coarse and fine b-splines represent the same identical spline curve.  Further inspection of~\eqref{eq:cpRefinement}  and subdivision matrix~\eqref{eq:subdivMatrix} reveals, that the fine control vertices are the weighted averages of coarse control vertices, with the \emph{columns} of the subdivision matrix representing the weights. There are two different sets of weights  corresponding to the two different types of columns of the subdivision matrix.  Based on the numbering scheme implied in Figure~\ref{fig:refinementCubicSpline} and the structure of the subdivision matrix, it is easy to deduce that one set of weights  applies to control vertices with even indices 
\begin{equation} \label{eq:evenIndices}
	\tilde {\vec x}_{2i} = \frac{1}{8} \vec x_{i-1} + \frac{3}{4} \vec x_i + \frac{1}{8} \vec x_{i+2} 
\end{equation}
and the other set of weights applies to control vertices with odd indices 
\begin{equation} \label{eq:oddIndices}
	\tilde {\vec x}_{2i+1} = \frac{1}{2} \vec x_i +  \frac{1}{2} \vec x_{i+1} \, .
\end{equation}
Hence, in terms of computer implementation,  for a given coarse polygon a corresponding fine polygon is obtained by first splitting each edge into two edges  and subsequently  computing control point coordinates according to~\eqref{eq:evenIndices} and~\eqref{eq:oddIndices}. In computer graphics literature the weights in~\eqref{eq:evenIndices} and~\eqref{eq:oddIndices} are usually given in form of stencils shown in Figure~\ref{fig:cubicMasks}.  

In subdivision schemes the foregoing described approach for obtaining a refined  control polygon is applied successively leading to finer and finer polygons. From the properties of the b-splines underlying~\eqref{eq:evenIndices} and~\eqref{eq:oddIndices},  it is clear that the control points converge to a cubic b-spline.

\subsection{Catmull-Clark subdivision surfaces}
\begin{figure*}
	\centering 
	 \null\hfill
	\subfloat[][Vertex] {
  		\includegraphics[scale=1]{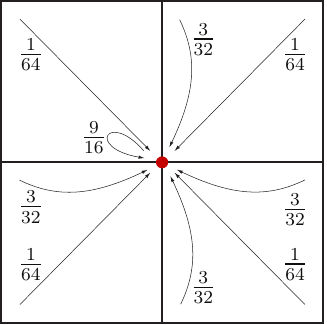}}
 	\hfill
  	\subfloat[][Face]{
		\includegraphics[scale=1]{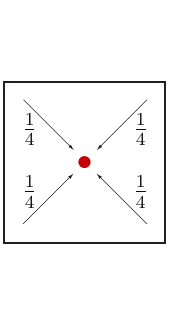}}
	\hfill	
 	\subfloat[][Edge]{
	  	\includegraphics[scale=1]{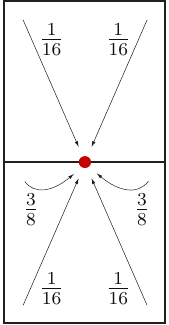}}
		 \hfill\null
  	\caption{Subdivision stencils for the Catmull-Clark scheme. Each of the stencils are used for computing the coordinates of vertices of the type indicated by red dot. }
  \label{fig:tensorProductMasks}
\end{figure*}

%
Two-dimensional b-splines can be generated as the tensor products of two one-dimensional b-splines. Similarly, the tensor product of two one-dimensional subdivision stencils yields the two-dimensional subdivision stencils shown in Figure~\ref{fig:tensorProductMasks}. The  three different stencils correspond to the three different type of vertices which occur during  the splitting of each face into four faces.  It can be easily verified that the weights given in Figure~\ref{fig:tensorProductMasks} are the tensor-products of the one-dimensional weights for cubic b-splines given in Figure~\ref{fig:cubicMasks}. Hence, successively refining a control mesh  and computing the control vertex coordinates with the stencils given in Figure~\ref{fig:tensorProductMasks} will lead in the limit to a cubic spline surface. 

\begin{figure*}
	\centering 
  		\includegraphics[scale=1]{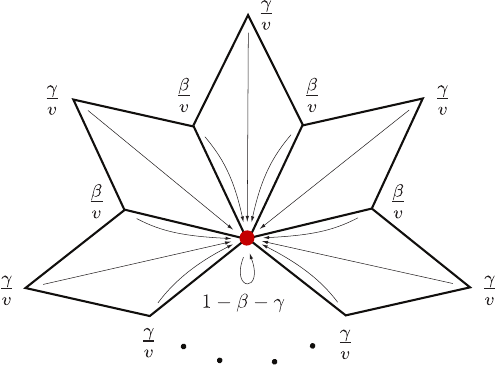}
	  	\caption{Subdivision stencil for an irregular vertex with valence $v$ with the weights  $\beta= \tfrac{3}{2v}$  and $\gamma= \tfrac{1}{4v}$, see~\cite{Catmull:1978aa}.}
  \label{fig:catmullClarkMask}
\end{figure*}

It is evident that the tensor-product stencils only apply to meshes in which each vertex within the domain is connected to four faces. The number of the faces connected to a vertex is referred to as the valence of that vertex and is denoted with~$v$. For the sake of brevity, we refer to~\cite{Zorin:2000aa} for the discussion of regularity of vertices on the boundaries and corners.  The domain vertices with a valence other than four are known as  extraordinary vertices or star-vertices. As originally proposed by Catmull and Clark~\cite{Catmull:1978aa},  the key idea in subdivision surfaces is to apply the modified stencil shown in Figure~\ref{fig:catmullClarkMask}  at the extraordinary vertices.   

To summarise, in each subdivision refinement step each face of the control mesh is split into four faces and the coordinates of the control vertices are computed with the weights given in Figures~\ref{fig:tensorProductMasks}  and~\ref{fig:catmullClarkMask} depending on the local connectivity structure of the vertex. There is mathematical theory which shows that the resulting surface is $C^2$ continuous almost everywhere except at the  extraordinary vertices where it is only $C^1$ continuous~\cite{Peters:2008aa}.

\begin{figure*}
	\centering 
	 \null\hfill
	 \subfloat[][Coarse control mesh ($\vec{x}^{\ell=0}$)] {
         \includegraphics[width=0.3\linewidth]{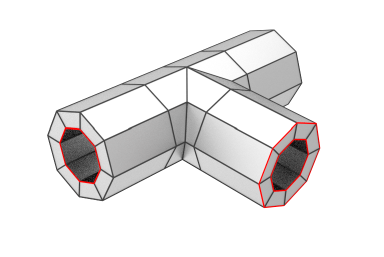}}
         \hfill
       \subfloat[][Once subdivided control mesh ($\vec{x}^{\ell=1}$)] {  
         \includegraphics[width=0.3\linewidth]{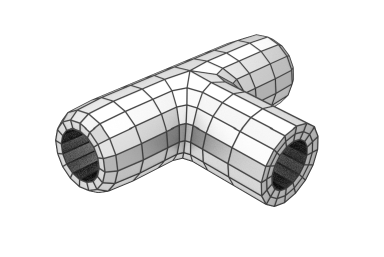}}
         \hfill
         \subfloat[][Limit surface ($\vec{x}^{\ell=\infty}$)] {
         \includegraphics[width=0.3\linewidth]{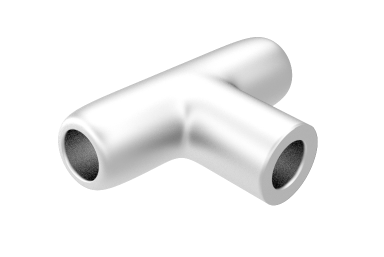}}
         \hfill\null
         \caption{Subdivision refinement of a T-junction geometry. On the coarse control mesh the edges in red are tagged as crease edges. The control mesh at the centre is obtained after one step of subdivision refinement. The geometry on the right is a rendering of the limit surface. Notice that the limit surface is not smooth on the crease edges.}
  \label{fig:subdivTjunction}
\end{figure*}

In addition to the stencils shown in Figures~\ref{fig:tensorProductMasks}  and~\ref{fig:catmullClarkMask} there are also extended subdivision stencils for vertices on edges, creases and corners, see, e.g.,~\cite{Ying:2001aa}. In this context, a crease is a line on the surface across which the surface is only $C^0$ continuous. As an illustrative example, Figure~\ref{fig:subdivTjunction} shows the subdivision refinement of a control mesh for a T-junction geometry  with extraordinary vertices and prescribed crease edges.

For later reference, we write the subdivision process as a linear mapping that maps a coarse control mesh at level $\ell$ to a finer control mesh at level $\ell+1$ 
\begin{equation} \label{eq:subdivisionAbstract}
	\vec x^{\ell+1} = \vec S \vec x^{\ell} \, ,
\end{equation}
where $\vec x^{\ell+1}$ and $ \vec x^{\ell}$ are two matrices containing  the coordinates of all the vertices at level $\ell$ and $\ell+1$, respectively. By definition the initial coarse control mesh has the level $\ell=0$.  The number of columns of $\vec x^{\ell+1}$ and $ \vec x^{\ell}$ is equal to the space dimension and their number of rows is equal to the number of all vertices in the mesh. The subdivision matrix $\vec S$ contains the weights given by the subdivision stencils and its dimension depends on the subdivision level $\ell$ considered.   Lastly, according to~\eqref{eq:subdivisionAbstract} the subdivision process can be interpreted as the chain of linear mappings for obtaining increasingly finer control meshes, i.e., 
\begin{equation}
	\begin{tikzpicture}[descr/.style={fill=white,inner sep=2.5pt}]
		\matrix (m) [matrix of math nodes, row sep=2.7em, column sep=2.4em, text height=1.5ex, text depth=0.25ex]	
		{ \vec{x}^0 & \vec{x}^1 &  \vec{x}^2  &   \cdots &  \vec{x}^{\ell-1}  &   \vec{x}^{\ell}  \\ };
		\path[->, font=\small]
		(m-1-1) edge node[above] {$ \vec S $} (m-1-2)
		(m-1-2) edge node[above] {$ \vec S $} (m-1-3)
		(m-1-3) edge node[above] {$ \vec S $} (m-1-4)
		(m-1-4) edge node[above] {$ \vec S $} (m-1-5)
		(m-1-5) edge node[above] {$ \vec S $} (m-1-6);
	\end{tikzpicture} \, .
\end{equation}
%

\subsection{Multiresolution surface editing \label{sec:multiresEdit}}
\begin{figure*}
	\centering 
         \includegraphics[width=0.9\textwidth]{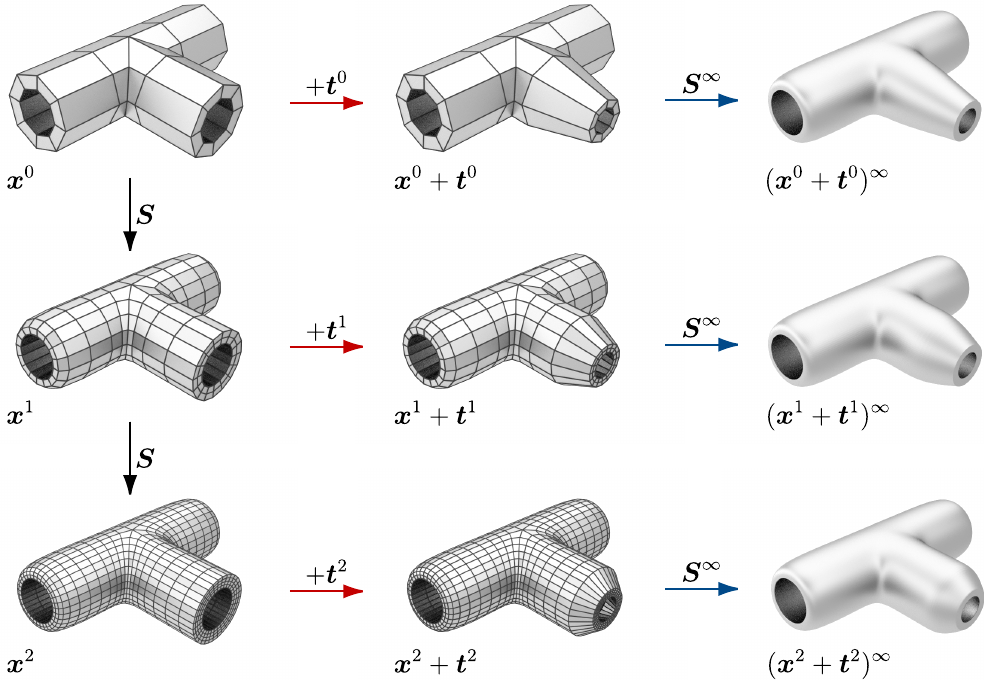}
         \caption{Multiresolution editing of the T-junction geometry  introduced in Figure~\ref{fig:subdivTjunction}.  The aim is to modify the shape and diameter of  the protruding tube. To achieve this the coordinates of selected vertices at either level  $\ell=0, 1 \text{ or } 2$  are modified as shown on the middle column, where the vectors $\vec t^0$, $\vec t^1$ and $\vec t^2$ represent the modifications at the respective level.  Notice the effect of the modification level on the limit surface (last column).}
  \label{fig:multiresEditing}
\end{figure*}

%
The sequence of control meshes generated during subdivision refinement readily lends itself for multiresolution editing of geometries~\cite{finkelstein1994multiresolution, Gortler:1995aa, zorin1997interactive, Lounsbery:1997aa}. As discussed,  Catmull-Clark subdivision surfaces are based on cubic b-splines. Hence, the support size of each subdivision basis function consists of a two-ring of faces, cf. Figure~\ref{fig:refinementCubicSpline}. With increasing refinement level the physical support size of basis functions becomes smaller. Accordingly, depending on the refinement level the editing of control vertex positions leads  to changes with different spatial extent on the limit surface. As an illustrative example the T-junction geometry introduced earlier is considered in Figure~\ref{fig:multiresEditing}.  In the middle column the coordinates of selected vertices at the levels $\ell=0$, $\ell=1$ and $\ell=2$ are modified. As can be seen, in the last column of pictures this leads  to changes in the limit surface in the vicinity of   the edited vertices and the  spatial extent of the changes is correlated with the refinement level. 

The subdivision surfaces by itself do not provide the possibility to simultaneously edit coarse and fine control meshes. For instance, after a fine control mesh is edited it is not possible anymore to edit a coarser level in order to apply large scale changes to the geometry. Simultaneous editing of different levels can be achieved with a wavelet-like  multiresolution decomposition of the control meshes, as will be discussed further below.

Before considering the multiresolution decomposition of  control meshes, we introduce  the coarsening of control meshes that were obtained with subdivision.  The linear  coarsening matrix $\vec R$ maps the given control points at level $\ell+1$ to the control points  at the coarser level $\ell$,
\begin{equation}	 \label{eq:costEquation}
	\vec x^\ell =  \vec{R} \vec x^{\ell+1}.
\end{equation}
The coarsening matrix $\vec R$ is not unique and different choices are possible. Essentially, the control mesh at level $\ell+1$ has more control vertices than the one at level $\ell$ and may contain more geometric information. In our implementation we obtain the coarsening matrix $\vec R$ from  a least squares fit of the subdivided  coarse control vertices $\vec S \vec x^\ell$  to the fine control vertices $\vec x^{\ell+1}$, i.e.,  
\begin{equation}
	\vec x^\ell = \underset{\vec y^\ell  }{\operatorname{argmin}}  \|   
  	\vec x^{\ell+1} - \vec{S} \vec y^\ell  \|^2 ,
\end{equation}
which leads to 
\begin{equation} \label{eq:lsqFit}
	\vec x^\ell = \vec R \vec  x^{\ell+1} \quad \text{ with } \vec{R}  = ({\vec{S}}^\trans \vec{S})^{-1} {\vec{S}}^\trans .
\end{equation}
By comparing with~\eqref{eq:subdivisionAbstract}  we can identify $\vec R$ as the pseudo-inverse of the subdivision matrix $\vec S$ so that
\begin{equation} 
	\vec x^\ell = \vec R \left (\vec S \vec x^\ell \right ) = \vec x^\ell .
\end{equation}
In words, subdivision refinement of a control mesh followed by coarsening (without editing) yields the original control mesh. Instead of least squares fitting, the coarsening matrix $\vec R$ can also be defined based on quasi-interpolation~\cite{deBoor1973,Litke:2001aa} or smoothing~\cite{zorin1997interactive}. On the other hand, coarsening by simply subsampling of the fine control mesh usually leads to artefacts in form of oscillations in the coarse control mesh. The proposed least squares fit approach is not very common in computer graphics because of the need for interactivity and fast processing times. Although the least squares matrix in~\eqref{eq:lsqFit} is sparse its solution cannot be found at interactive rates.  Notice also that~\eqref{eq:lsqFit}  represents a system of equations with $d$ right hand side vectors for each of the $d$ coordinate directions.  

\begin{figure*}
	\centering 
	 \null\hfill
	 \subfloat[][Fine control mesh ($\vec{x}^{\ell=3}$)] {
         \includegraphics[width=0.32\linewidth]{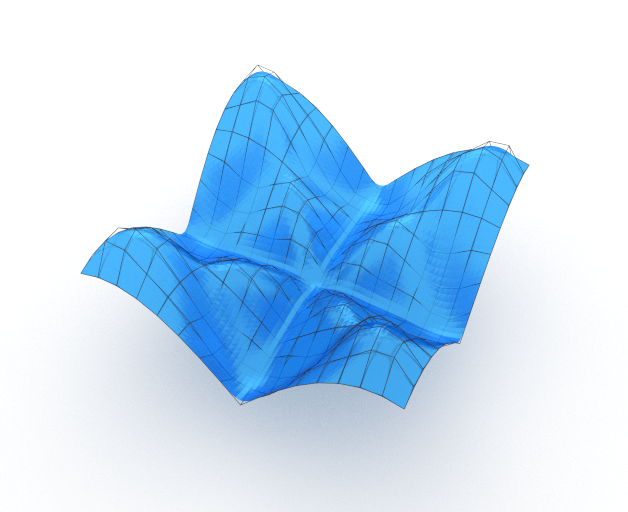}}
         \hfill
       \subfloat[][Twice coarsened control mesh ($\vec{x}^{\ell=1}$)] {  
         \includegraphics[width=0.32\linewidth]{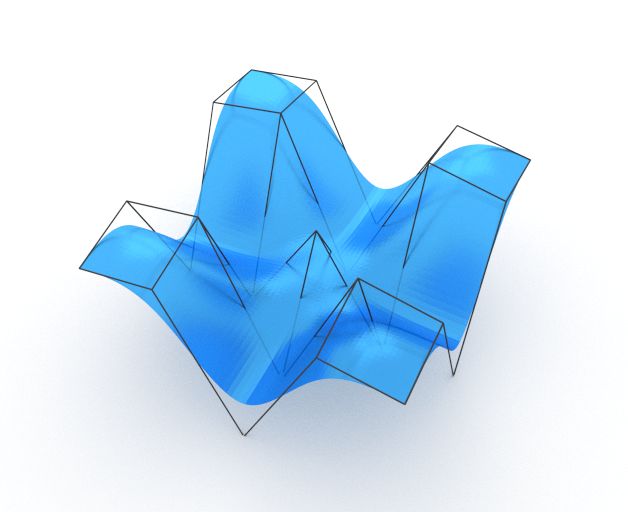}}
         \hfill
         \subfloat[][Coarse control mesh ($\vec{x}^{\ell=0}$)] {
         \includegraphics[width=0.32\linewidth]{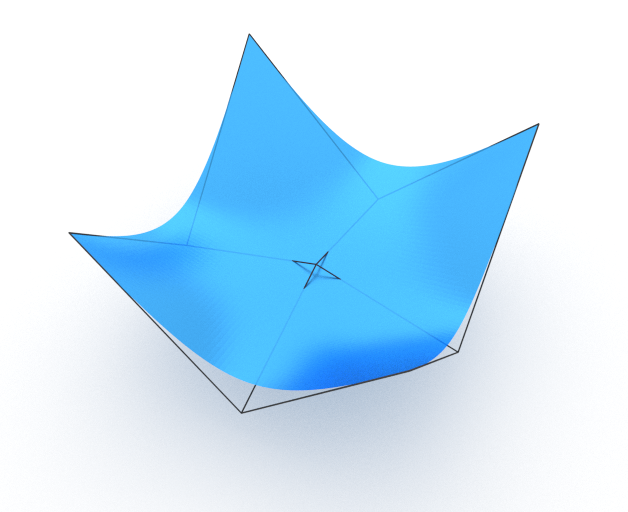}}
         \hfill\null
         \caption{Subdivision coarsening of a fine control mesh. Successive application of the coarsening matrix $\vec R$ leads to increasingly coarser control meshes. }
  \label{fig:subdivisionCoarsening}
\end{figure*}

Similar to  subdivision refinement the coarsening matrix can be successively applied in order to obtain coarser representations of the geometry, i.e., 
\begin{equation}
	\begin{tikzpicture}[descr/.style={fill=white,inner sep=2.5pt}]
		\matrix (m) [matrix of math nodes, row sep=2.7em, column sep=2.4em, text height=1.5ex, text depth=0.25ex]	
		{ \vec{x}^0 & \vec{x}^1 &  \vec{x}^2  &   \cdots &  \vec{x}^{\ell-1}  &   \vec{x}^{\ell}  \\ };
		\path[<-, font=\small]
		(m-1-1) edge node[above] {$ \vec R $} (m-1-2)
		(m-1-2) edge node[above] {$ \vec R $} (m-1-3)	
		(m-1-3) edge node[above] {$ \vec R $} (m-1-4)
		(m-1-4) edge node[above] {$ \vec R $} (m-1-5)
		(m-1-5) edge node[above] {$ \vec R $} (m-1-6);
	\end{tikzpicture} \, .
\end{equation}
The dimension of the matrix $\vec R$ depends on the considered level~$\ell$. As mentioned during the coarsening process each of the coordinate directions are considered individually.  Figure~\ref{fig:subdivisionCoarsening} shows the coarsening of a subdivision surface with the described approach. From the shown limit surfaces it is visible that the coarsening process leads to a smoothing of the geometry; and the overall geometry is faithfully represented by the coarser representations.  For this reason, the coarsening process is sometimes also referred to as the smoothing process. 

With the introduced subdivision and coarsening operations it is now possible to  devise a wavelet-like multiresolution decomposition of a control mesh and the associated geometry. The aim of this decomposition is to enable the simultaneous editing of different levels of the subdivision surface. This is achieved by    storing a coarse control mesh and the  differences between successive control meshes called \emph{details}. Subsequently it is possible to first edit any control mesh level and then to add the stored details as necessary.  For a given fine control mesh at level $\ell$ the detail vector  $\vec d^{\ell-1}$ is computed using the subdivision and coarsening processes as follows 
\begin{equation} \label{eq:analysis}
	\vec d^{\ell-1} = \vec x^\ell - \vec S \left ( \vec R  \vec x^\ell \right ) = (\vec I - \vec S \vec R ) \vec x^\ell \, .
\end{equation}
In turn, when a control mesh at level ${\ell-1}$  and the detail $\vec d^{\ell-1}$ are given the geometry at level $\ell$ can be recovered according to~\eqref{eq:analysis} and \eqref{eq:lsqFit} with 
\begin{equation} \label{eq:synthesis}
	\vec x^\ell =  \vec S \vec x^{\ell-1} + \vec d^{\ell-1}  \, .
\end{equation}
The global detail vector $\vec d^\ell$ is composed of local vertex detail vectors, which can be conveniently stored at the vertices. In our actual implementation, we express each local detail vector in a local tangential coordinate system at its vertex. As known in computer graphics, this is necessary so that in~\eqref{eq:synthesis} any modifications to the geometry $\vec x^{\ell-1}$ should have an intuitive effect on the details contained in $\vec d^\ell$, see, e.g.,~\cite{finkelstein1994multiresolution, zorin1997interactive, Sorkine:2004aa}.

Finally, the two-level decomposition given by \eqref{eq:analysis} and \eqref{eq:synthesis} can be successively applied leading to a wavelet-like multiresolution decomposition of the surface. The process of obtaining the details for a given fine geometry is referred to as the \emph{analysis step} 

\begin{equation} \label{eq:waveletAnalDiag}
	\begin{tikzpicture}[descr/.style={fill=white,inner sep=2.5pt}]
		\matrix (m) [matrix of math nodes, row sep=2.7em, column sep=2.4em, text height=1.5ex, text depth=0.25ex]	
		{ \vec{x}^0 & \vec{x}^1 &  \vec{x}^2  &   \cdots &  \vec{x}^{\ell-1}  &   \vec{x}^{\ell}  \\
		 \vec{d}^0 & \vec{d}^1 & \cdots & \vec{d}^{\ell-2}  & \vec{d}^{\ell-1} &   \\ };
		\path[->, font=\small]
		(m-1-2) edge node[above] {$ \vec R $} (m-1-1)
		(m-1-2) edge node[below,sloped] {$ \vec{I} - \vec{SR} $} (m-2-1)
		(m-1-3) edge node[above] {$ \vec R $} (m-1-2)
		(m-1-3) edge node[below,sloped] {$ \vec{I} - \vec{SR} $} (m-2-2)
		(m-1-4) edge node[above] {$ \vec R $} (m-1-3)
		(m-1-4) edge node[below,sloped] {$ \vec{I} - \vec{SR} $} (m-2-3)
		(m-1-5) edge node[above] {$ \vec R $} (m-1-4)
		(m-1-5) edge node[below,sloped] {$ \vec{I} - \vec{SR} $} (m-2-4)
		(m-1-6) edge node[above] {$ \vec R $} (m-1-5)
		(m-1-6) edge node[below,sloped] {$ \vec{I} - \vec{SR} $} (m-2-5);
	\end{tikzpicture}
\end{equation}
The corresponding \emph{synthesis step} takes the form 
\begin{equation} \label{eq:waveletSyntDiag}
\begin{tikzpicture}[descr/.style={fill=white,inner sep=2.5pt}]
\matrix (m) [matrix of math nodes, row sep=2.7em, column sep=2.4em, text height=1.5ex, text depth=0.25ex]	
	{ \vec{x}^0 & \vec{x}^1 &  \vec{x}^2  &   \cdots &  \vec{x}^{\ell-1}  &   \vec{x}^{\ell}  \\
	 \vec{d}^0 & \vec{d}^1 & \cdots & \vec{d}^{\ell-2}  & \vec{d}^{\ell-1} &   \\ };
	\path[->, font=\small]
	(m-1-1) edge node[above] {$ \vec S $} (m-1-2)
	(m-2-1) edge (m-1-2)
	(m-1-2) edge node[above] {$ \vec S $} (m-1-3)
	(m-2-2) edge  (m-1-3)
	(m-1-3) edge node[above] {$ \vec S $} (m-1-4)
	(m-2-3) edge (m-1-4)
	(m-1-4) edge node[above] {$ \vec S $} (m-1-5)
	(m-2-4) edge (m-1-5)
	(m-1-5) edge node[above] {$ \vec S $} (m-1-6)
	(m-2-5) edge (m-1-6);
\end{tikzpicture}
\end{equation}
For an efficient implementation of the analysis and synthesis steps and the related data structures we refer to Zorin et al.~\cite{Zorin:2000aa}.

\begin{figure*}
	\centering 
         \includegraphics[width=0.85\linewidth]{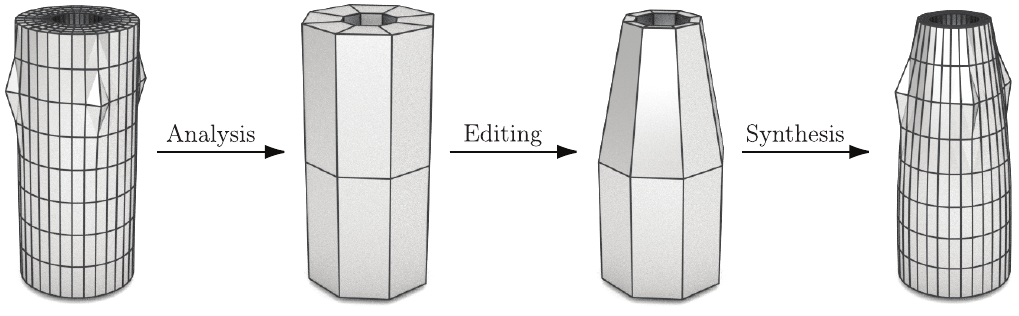}
         \caption{Multiresolution editing of a cylindrical component. The fine resolution mesh shown on the left at level $\ell=2$ is given.  In the analysis step the geometry is decomposed into the  coarse level $\ell=0$ (second from the left) and details $\vec d^0$ and $\vec d^1$, cf.~\eqref{eq:waveletAnalDiag}.  As shown on the third from the left, the coarse level can be edited irrespective of the details. In the subsequent  synthesis step~\eqref{eq:waveletSyntDiag} adding the precomputed details  to the coarse level yields the control mesh shown on the right.  }
  \label{fig:cylinderEditing}
\end{figure*}

A typical workflow during multiresolution editing of  a cylindrical component  with few small bumps is  shown in Figure~\ref{fig:cylinderEditing}. First the geometry is created starting from a coarse control mesh and adding the bumps on the two times subdivided control mesh. This is achieved by displacing few selected vertices on the fine control mesh and yields the leftmost picture Figure~\ref{fig:cylinderEditing}. After this step, in order to apply large scale changes it is necessary to employ a multiresolution decomposition of the geometry. More specifically, the control mesh $\vec x^{\ell=2}$ is decomposed into the detail vectors~$\vec d^0$ and $\vec d^1$, and the original control mesh~$\vec{x}^0$.  For this specific geometry the detail vector $\vec d^0$ is zero.  After this decomposition it is possible to  change the original control mesh into, for instance, a cone and subsequently to subdivide and to automatically add the stored details leading to the shown cone geometry with bumps. 

\subsection{Multiresolution shape optimisation \label{sec:multiresOpt}}
%
The introduced subdivision multiresolution editing technique enables the use of two different resolutions of the same geometry for optimisation and analysis. The two resolutions correspond to different  refinement levels in a multiresolution  hierarchy. In shape optimisation it is usually necessary  to use  a  coarse control mesh  for geometry updating and a relatively fine control mesh for analysis.  As is known, unwanted geometry oscillations may appear when the  analysis and geometry  representations have similar resolutions \cite{Braibant:1984aa,le2011gradient,bletzinger2014consistent}. These geometry oscillations are usually a numerical artefact or a result from the ill-posedness of the considered optimisation problem. Moreover, in practical applications it might be desirable to optimise only a very coarse representation   out of aesthetic or manufacturability reasons.

The Algorithm~\ref{algo:multiresTopDown} describes the proposed multiresolution shape optimisation approach. The fixed computational level $\ell_c$ and the maximum optimisation level $\ell_{o,max}$ are prescribed by the user. The level $\ell_c$ has to be large enough such that the numerical solution is accurate enough for practical purposes. The maximum optimisation level has to be chosen $ \ell_{o,max} \le  \ell_c$ and determines the smallest geometric feature size on the optimised geometry. The input control mesh $\vec x^{\ell_{inp}}$ can be a coarse mesh with $\ell_{inp}=0$ or an already edited fine multiresolution mesh with $\ell_{inp}>0$.  For control meshes with $\ell_{inp}>0$ first a multiresolution decomposition as indicated in~\eqref{eq:waveletAnalDiag} is performed. Throughout the algorithm the optimisation control mesh and its level are denoted with  $\vec x^{\ell_o}$ and $\ell_o$, respectively. For the immersed finite element analysis the geometry corresponding to a  control mesh $\vec x^{\ell_c}$  with $\ell_c \ge \ell_o$ is used. The analysis level $\ell_c$ is usually fixed and the control mesh has elements of similar  size like the cells  of the immersed finite element grid. In order to obtain the analysis control mesh $\vec x^{\ell_c}$ from the optimisation control mesh $\vec x^{\ell_o}$ we use the introduced multiresolution refinement technique.  The immersed finite element analysis yields the cost function $ J (\vec x^{\ell_c}, \vec u(\vec x^{\ell_c})) $ and the shape kernel $g(\vec x^{\ell_c})$, see~\eqref{eq:shapeDerivFinal}.  To compute the shape gradient at the vertices, the surface normal vector $\vec n (\vec x^{\ell_c})$ is required, which is easily  computed  with the known subdivision stencils for tangent vectors, see e.g.~\cite{Zorin:2000aa}. Vertex-wise multiplication of the shape kernel with the normal vector yields the shape gradient vector $\vec g^{\ell_c}$, which is subsequently projected to the optimisation level vector $\vec g^{\ell_o}$  by successive coarsening with $\vec R$.  For the sake of brevity, in Algorithm~\ref{algo:multiresTopDown} the geometry is updated with a steepest descent technique and no additional constraints are present. In applications it is common to have additional constraints, such as perimeter, area or volume constraints. In our  actual implementation  we optimise the constrained discrete problem with the method of moving asymptotes (MMA) proposed by Svanberg~\cite{Svanberg:1987, Svanberg:2002aa} using the implementation in  the NLopt library~\cite{nloptPackage}. Finally, note that in Algorithm~\ref{algo:multiresTopDown} the optimisation level $\ell_o$ is not fixed, it is incremented after a minimum is reached and until a user prescribed maximum optimisation level is reached $\ell_o \le \ell_{o,max}$.   

\begin{algorithm}[H]
\caption{Multiresolution shape optimisation}
\label{algo:multiresTopDown}
\begin{algorithmic}[1]
\LCOMMENT{choose maximum optimisation level~$\ell_{o,max}$ and computational level~$\ell_c$ }
\LCOMMENT{read input control mesh~$\vec x^{\ell_{inp}}$ }
\IF{$\ell_{inp}=0$} 
	\LCOMMENT{initialise all detail vectors ($\vec d^{\ell} = \vec 0$) }
\ELSIF {$\ell_{inp} >  0$} 
 	\FOR{$ \ell \gets \ell_{inp} \textrm{ to } 0$}
 		\STATE $\vec d^{\ell} = \vec x^{\ell} - \vec S (\vec R \vec x^{\ell})$
 	\ENDFOR 
\ENDIF
\LCOMMENT{Initialise optimisation level}
\STATE $\ell_o =  0$  
\LCOMMENT{Initialise cost function}
\STATE $J =  \infty$  
\LCOMMENT{iterate over optimisation levels}
\WHILE {$\ell_o \le \ell_{o,max}$}
\LCOMMENT{update vertex coordinates $\vec x^{\ell_o}$ while the cost function decreases}
\REPEAT 
\LCOMMENT{subdivide optimisation level $\ell_o$ up to analysis level $\ell_c$}
 \FOR{$ \ell \gets \ell_o \textrm{ to } \ell_c$}
 	\STATE $\vec{x}^{\ell} \gets \vec S \vec x^{\ell} + \vec d^{\ell} $
 \ENDFOR 
\LCOMMENT{compute cost function $J = J (\vec x^{\ell_c}, \vec u(\vec x^{\ell_c})) $ and shape kernel field $ g(\vec x)  = 
g (\vec u (\vec x^{\ell_c}))$}

\LCOMMENT{compute maximum ascent direction at the vertices $\vec x^{\ell_c}$ with the outer normals $\vec n ( {\vec x^{\ell_c}})$ }

\STATE $\vec g^{\ell_c} =   g (\vec x^{\ell_c})  \vec n ( {\vec x^{\ell_c}})$

\LCOMMENT{project shape derivative to optimisation level}

\FOR{$ \ell \gets \ell_c \textrm{ to } \ell_o$}
 	\STATE $\vec{g}^{\ell} \gets \vec R \vec g^{\ell}  $
 \ENDFOR 
 
\LCOMMENT{update vertex coordinates of the optimisation level}
\STATE $\vec x^{\ell_o}  \leftarrow  ( \vec x ^{\ell_o} - \alpha  \vec g^{\ell_o}) \quad
\text{ with } \alpha \ge 0$
\UNTIL{{$ ( J_{\text{previous}} -J ) <  \text{tolerance}$}}

\LCOMMENT{increment optimisation level}
\STATE $\ell_o \gets (\ell_o+1)$
\STATE $\vec x^{\ell_o} \gets  \vec S \vec x^{\ell_o} + \vec d^{\ell_o}$
\ENDWHILE
\end{algorithmic}
\end{algorithm}

\section{Examples \label{sec:examples}}
In this section, we present several examples to demonstrate the robustness  and versatility of the proposed multiresolution shape optimisation technique. In all the examples the domain is described either by a cubic b-spline curve (in 2D) or a Catmull-Clark subdivision surface (in 3D). The objective of the optimisation is to minimise the structural compliance~\eqref{eq:compliance}, which is equivalent to maximising the structural stiffness. The corresponding adjoint problem~\eqref{eq:adjointStrong} has (up to the sign) the same right hand side as the primal problem~ \eqref{eq:elasticityStrong}. Hence, it is sufficient to consider only the primal problem, which is solved with the immersed finite element technique using quadratic b-spline basis functions. The resulting smooth stress field in combination with unique  boundary normals at the vertices of the control mesh leads to smooth shape gradients. The optimised boundary curves or surfaces have in general no corners or sharp edges so that there is always a unique normal.  

Initially, we consider in Section~\ref{sec:shapeOptEx} only shape optimisation examples. Subsequently, in Section~\ref{sec:shapeTopOptEx} in addition to the shape also the topology of the domain is optimised.
 To this end, we make use of the topology derivative, see e.g.~\cite{eschenauer2001topology, novotny2003topological}, to introduce new holes in the domain. In our current implementation the merging or removing of holes is  not considered.
In all the two-dimensional examples we use a plane strain formulation unless otherwise indicated. 

\subsection{Shape optimisation \label{sec:shapeOptEx}}

\subsubsection{Simply supported plate with a hole}
\begin{figure*}
  \centering 
 	\subfloat[][Problem description.] {
	\includegraphics[scale=0.9]{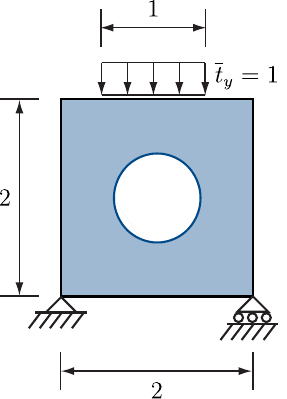}} \hspace{0.05\textwidth}
	\subfloat[][Optimised hole shapes for cases C1, C2 and C3.] {
    \includegraphics[scale=0.9]{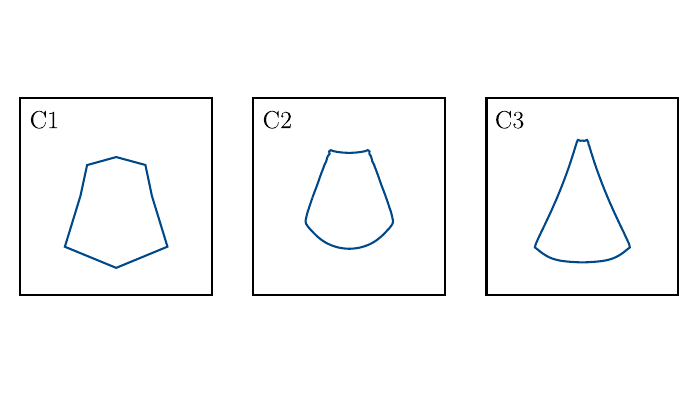}}
   \caption{Simply supported plate with a hole.  \label{fig:optMotivationA}}
\end{figure*}

%
This introductory example aims to highlight the advantages of multiresolution optimisation over classical approaches using only one or two representation levels. The problem consists of a square plate with an edge length $L=2$ and a circular hole with diameter $D=1$, see Figure~\ref{fig:optMotivationA}. The plate is loaded with a line load of length of $1$.  The Young's modulus and Poisson's ratio of the plate are $E=100$ and $\nu = 0.4$, respectively. 

During optimisation the shape of the hole is to be modified so that the structural compliance of the plate is minimised. The area of the hole is constrained with  
\begin{equation} \label{eq:areaConst}
	A (\vec x^{\ell_c}) \ge \frac{\pi}{4} \, ,
\end{equation}
where $\vec x^{\ell_c}$ are the vertex coordinates of the analysis control mesh at level $\ell_c$.  This constraint is necessary since the stiffest plate is the one with a zero hole diameter. The area of the hole is computed by integrating over its boundary
\begin{equation} \label{eq:bloxArea}
	A(\vec x^{\ell_c}) = \frac{1}{2} \int_{\Gamma} \vec x^{\ell_c} \cdot \vec n (\vec x^{\ell_c}) \dif \Gamma  \, ,
\end{equation}
where $\vec n (\vec x^{\ell_c}) $ is the normal to the boundary curve $\Gamma$. Recall that  we represent the boundary curve with cubic b-splines. In computations, \eqref{eq:bloxArea} is evaluated using one quadrature point per element of the control polygon. During the optimisation also the gradient of the area constraint is required, which is computed by differentiating the discretised version of~\eqref{eq:bloxArea} with respect to the vertex positions of the analysis control mesh $\vec x^{\ell_c}$.  Alternatively, it would be possible to use analytic shape derivatives equivalent to~\eqref{eq:reynolds}. The vertex-wise gradient of the area constraint on the analysis level $\ell_c$ is projected to the optimisation level $\ell_o$ in the same way as the shape gradient vector $\vec g^{\ell_c}$. The detail vectors created during the decomposition~\eqref{eq:waveletAnalDiag} are discarded. 

Initially, at level $\ell=0$ the hole is represented with a cubic spline with $8$ control points. The immersed finite element grid has $100 \times 100$ cells of uniform size. Three cases referred to as C1, C2 and C3 with different geometry and analysis resolutions are studied: 
\begin{itemize}
	 \setlength{\itemsep}{1pt}
  	\setlength{\parskip}{0pt}
  	\setlength{\parsep}{0pt}
	\item [-] In C1 only one level with $\ell_o = \ell_c=0$ is used for analysis and optimisation.
	\item [-] In C2 a four times subdivided control mesh at refinement level $\ell_o = \ell_c=4$ is used for analysis and optimisation.
	\item [-] In C3 the optimisation level starts with $\ell_o=0$ and increases until $\ell_o = \ell_c = 4$ is reached. Throughout the computations the analysis level is fixed to $\ell_c=4$
\end{itemize}
In case C1 the control mesh that is visible by the immersed finite element grid contains $8$ elements and in cases C2 and C3 it contains $128$ elements. It is clear that in case C1 the hole geometry is poorly resolved on the immersed finite element grid.

\begin{figure*}
	\centering
	\includegraphics[scale=1]{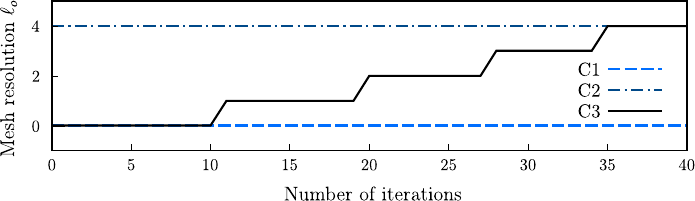}
	\caption{Simply supported plate with a hole. Variation of the optimisation level $\ell_o$ over the number of optimisation iterations. \label{fig:optMotivationB}}
\end{figure*}

\begin{figure*}
	\centering 
	\includegraphics[scale=1]{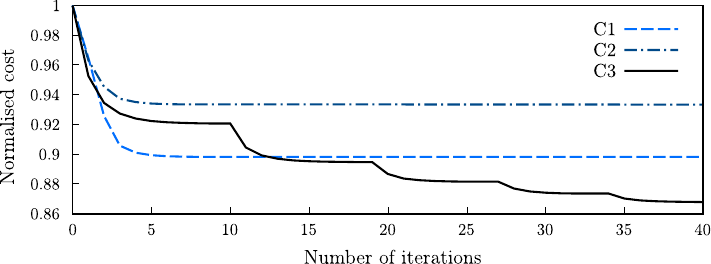}
	\caption{Simply supported plate with a hole. Reduction of the normalised cost over the number of optimisation iterations. The initial cost for case C1 is $0.073$ and for cases C2 and C3 is $0.065$. \label{fig:optMotivationC}} 
\end{figure*}

In Figure~\ref{fig:optMotivationA} the optimised final hole shapes for the three cases are shown. In particular, the difference in optimal shapes for cases C2 and C3 which use the same analysis level $\ell_c=4$ is striking. The case C1 is different from the other two cases because of the mentioned inadequately coarse analysis control mesh with $\ell_c =0$. As indicated in Figure~\ref{fig:optMotivationB}, during optimisation only for case C3 the optimisation level $\ell_o$ is successively increased.  The optimisation level is always incremented when a minimum is reached, cf. Algorithm~\ref{algo:multiresTopDown}.   For the three cases the reduction of the relative cost function over the number of iterations is shown in Figure~\ref{fig:optMotivationC}.  The case C2 with fixed fine resolution achieves the smallest cost reduction while the case C3 with multiresolution achieves the largest cost reduction. The strong dependence of the optimisation results on geometry parameterisation is well known in structural optimisation and is often associated with the non-convexity of the considered optimisation problem. We conjecture that by initially using a coarse control mesh for optimisation  the possible number of local minima is significantly reduced which reduces the possibility of landing in a non-optimal local minimum. It appears that in case C2 the optimisation problem is caught in a local minimum which is significantly higher than the global minimum. 

\subsubsection{Optimal hole shapes in a two-dimensional domain \label{sec:holesInPlate2D}}
\begin{figure}
  \centering 
   \includegraphics[scale=1]{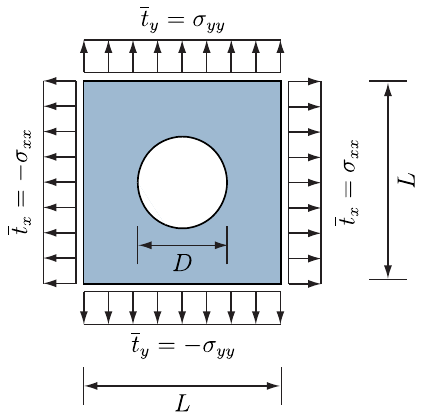}
  \caption{Optimal hole shapes in a two-dimensional domain. Problem description. \label{fig:cavity2Dsystem}} 
\end{figure}

%
In this prototypical example we study the optimal hole shapes in an elastic domain subjected to bi-axial stress.  The problem setup is shown in Figure~\ref{fig:cavity2Dsystem}. As in the previous example the area of the hole is constrained to be constant during optimisation. The Young's modulus and Poisson's ratio are chosen with $E=100$ and $\nu=0.4$, respectively.  According to analytical results for infinite plates the optimal hole shape depends on the ratio and sign of the far-field stress~\cite{Cherkaev:1998aa}. When the two components of the far-field stress are of the same sign the optimal hole shape is an ellipse with an aspect ratio $r_x/r_y$ equal to the bi-axial stress ratio  $\sigma_{xx}/\sigma_{yy}$. In contrast, when the two far-field stress components are of opposite sign the optimal hole shape is a quadrilateral with smooth corners. The case with far-field stresses of the same sign has been widely studied in literature, see, e.g.,~\cite{kristensen1976optimum, Norato:2004aa, cervera2005evolutionary}.

In our computations the initial hole geometry at level $\ell=0$ is modelled with a cubic spline with $8$ control vertices. The position of the vertices is chosen such that the resulting spline curve represents approximately a circle with diameter $D=1$.  The three times subdivided control mesh with $64$ vertices serves as the computational mesh for describing the hole geometry on the immersed finite element grid.  The optimisation starts with $\ell_o=0$ and is incremented until $\ell_o = \ell_c = 3$ is reached, cf. Algorithm~\ref{algo:multiresTopDown}. 

\begin{figure}
  \centering 
  \includegraphics[scale=1]{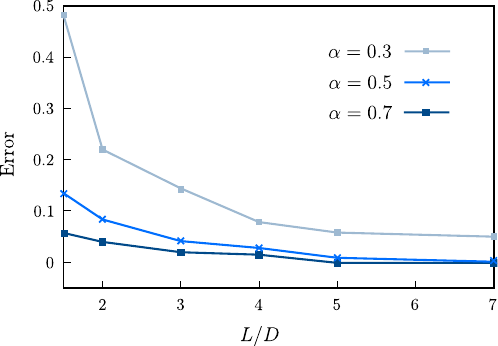}
  \caption{Optimal hole shapes in a two-dimensional domain. The relative error in the computationally obtained aspect ratio of the ellipse for different stress ratios $\alpha$ and plate sizes $L/D$, with the hole diameter $D=1$. The error is defined as the difference between the computationally  and analytically obtained aspect ratios.  \label{fig:cavity2Dconv}}
\end{figure}

In a first set of computations we quantify the effect of computing with a finite size domain as opposed to an infinite domain underlying the analytical results. To this end, the length to diameter ratio $L/D$ is varied between $1.5 \leq L/D \leq 7$  while the hole diameter is fixed with $D=1$. In all computations the element size on the immersed finite element grid is kept fixed with $1/25$. We compute the optimised hole shapes for three different bi-axial stress ratios $\alpha = \sigma_{xx}/\sigma_{yy} \in  \{0.3, 0.5, 0.7 \}$ and seven different length to diameter ratios $L/D \in \{1.5, 2, 3, 4, 5, 7\}$. As mentioned, for bi-axial  stress components with the same sign the analytically obtained optimal hole shape is an ellipse with an aspect ratio equal to the stress ratio~$\alpha$. Figure~\ref{fig:cavity2Dconv} shows the error in the computationally obtained ellipse aspect ratios for different $\alpha$ and $L/D$. The error is due to the finite size of the domain and the discretisation errors. 
 As can be seen in Figure~\ref{fig:cavity2Dconv} for stress ratios $\alpha \in \{0.5, 0.7\}$ the computationally obtained aspect ratio converges to the analytic result  for sufficiently large domains. However, for $\alpha=0.3$ the obtained aspect ratio does not converge to the analytic result. This is due to the large discretisation errors close to the two apexes of the relatively tall ellipse. This error could be reduced by increasing the computational level $\ell_c$ and decreasing the cell size of the immersed finite element grid.

In the following set of computations the domain size  and the initial hole diameter are chosen with $L=4$ and $D=1$. The cell size of the immersed finite element grid is chosen with $1/100$. According to Figure~\ref{fig:cavity2Dconv} this setup appears to provide sufficient accuracy while keeping the computation time manageable. The geometric description of the hole remains the same as in the previous set of computations. With the described setup we compute the optimised hole shapes for ten different stress ratios 
\[ \alpha \in\{-1.0, -0.7, -0.5, -0.3, -0.1, 0.1, 0.3, 0.5, 0.7, 1.0\} \, . \] 
Cherkaev et al.~\cite{Cherkaev:1998aa}  has shown that for negative stress ratios having several holes gives a lower compliance than having one single hole. Therefore, it is necessary to prevent  the appearance of multiple holes which cannot be obtained with shape optimisation starting with a single hole. To regularise the optimisation problem we add to the compliance cost function~\eqref{eq:compliance} an integral penalising perimeter change, i.e., 
\begin{equation} \label{eq:compliancePenalty}
	J (\Omega, \vec u ) =  \int_{\Omega} \vec \sigma (\vec u ) : \vec \epsilon (\vec u) \,  \dif \Omega   + \rho_L \left ( \int_\Gamma d\Gamma \right )^2  \, .
\end{equation}
where $\rho_L $ is a prescribed penalty parameter. The parameter $\rho_L$ controls how much regularisation is applied to penalise  the formation of multiple holes.  Without this regularisation the solution of the discretised optimisation problem with the method of moving asymptotes (MMA) does not converge. In computations we integrate the  second term in~\eqref{eq:compliancePenalty} numerically using the control mesh on the computational level~$\ell_c$. The required gradient information is obtained by differentiating the resulting discrete equations with respect to the vertex positions.  The obtained gradient vector is added to the shape gradient vector~$\vec{g}^{\ell_c}$.  

\begin{figure*}
  \centering 
  {\includegraphics[scale=1]{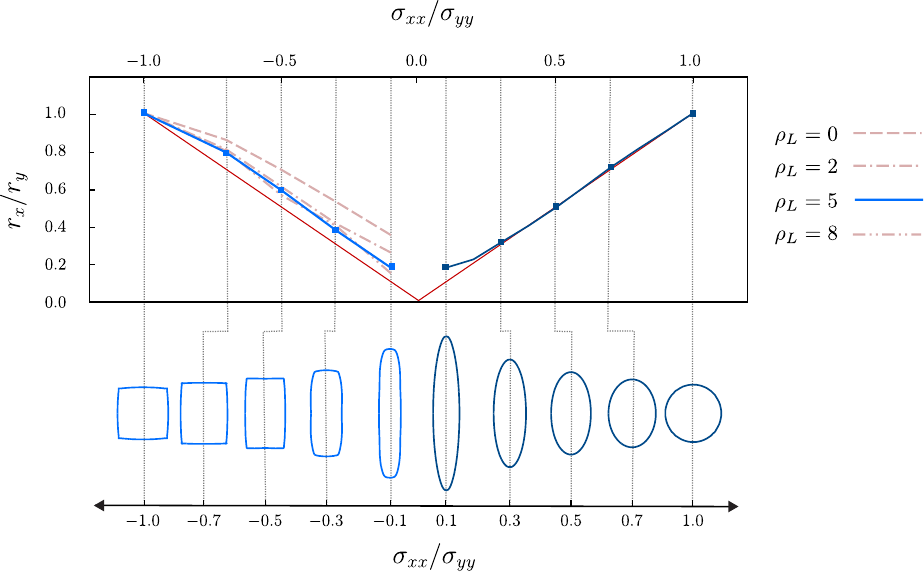}}
  \caption{Optimal hole shapes in a two-dimensional domain.  Computationally obtained hole shapes and their aspect ratios $r_x / r_y$ for different stress ratios $\alpha=\sigma_{xx} / \sigma_{yy}$.  The analytically obtained aspect ratio is shown in red. The multiple curves for  $\sigma_x / \sigma_y<0$ are computed using different  penalty values $\rho_L$, cf.~\eqref{eq:compliancePenalty}.  No penalty is applied when  $\sigma_x / \sigma_y>0$.}
  \label{fig:cavityShapesB}
\end{figure*}

With cost function~\eqref{eq:compliancePenalty} it is now possible to compute the optimal hole shapes for positive as well as negative stress ratios. The obtained hole shapes  and aspect ratios are shown in Figure~\ref{fig:cavityShapesB}. For positive stress ratios the hole is of elliptical shape and for negative ratios it is a smoothed quadrilateral. The area of all the holes is equal due to the prescribed area constraint. As can be seen in Figure~\ref{fig:cavityShapesB},  the obtained aspect ratios are in very good agreement with the analytical result indicated by the solid red line. The discrepancy for very small stress ratios is due to the finite size of the domain and the appearance of tall holes with crack-like stress concentrations requiring a finer discretisation. Finally,  the effect of the penalty parameter $\rho_L$ on the obtained aspect ratio is relatively mild. 
%

%
\subsubsection{Optimal hole shapes in a three-dimensional domain \label{sec:holesInPlate3D}}
%
The considered computational domain is a cube with a side length of $L=4$ and is discretised with cells of size $1/20$.   The initial geometry of the hole is a sphere with diameter $D=1$ and is at level $\ell=0$  approximated with a control mesh of $26$ nodes. The optimisation level starts with $\ell_o=0$ and is incremented until $\ell_o = \ell_c = 3$ is reached. The volume of the hole is constrained to remain constant during optimisation. This is achieved by computing the volume and its gradient with the three-dimensional extension of~\eqref{eq:bloxArea}.

\begin{figure}
  \centering 
  {\includegraphics[scale=0.7]{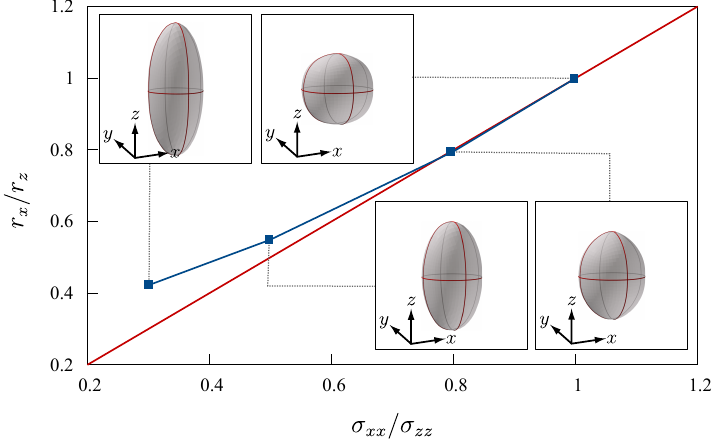}}
  \caption{Optimal hole shapes in a three-dimensional domain. Computationally obtained ellipsoidal hole shapes and their aspect ratios $r_x/r_z$ for different stress ratios $\sigma_{xx}/\sigma_{zz}$. Note that $\sigma_{xx} = \sigma_{yy}$ and $r_x = r_y$ (up to discretisation errors). The analytically obtained aspect ratio is shown in red.  \label{fig:optCavityShapesAndRatio3D}}
\end{figure}

First we choose the two stress components  $\sigma_{xx} =\sigma_{yy}$ as equal and only modify the $\sigma_{zz}$ stress component such that 
\[ 
\sigma_{xx}/\sigma_{zz}\in\{0.3,0.5,0.8,1.0\}. \]
According to analytical results the optimised hole shapes are ellipsoids with semi-axis radius $r_x = r_y$ and have the aspect ratio $r_x/r_z = \sigma_{xx} / \sigma_{zz}$.  Figure~\ref{fig:optCavityShapesAndRatio3D}  shows the computationally obtained ellipsoidal hole shapes and their aspect ratios. The computational and analytical results agree very well especially for large stress ratios.   For smaller stress ratios the discrepancy is due to the discretisation errors in resolving the more pronounced stress concentrations caused by taller holes.  

Next, we keep the stress component $\sigma_{zz}$  fixed and independently vary the two stress components $\sigma_{xx}$ and $\sigma_{yy}$. To reduce the effect of  domain size on the computationally  obtained hole shapes we chose the domain dimensions in dependence of the stress ratios $\alpha_x = \sigma_{xx}/\sigma_{zz}$ and $\alpha_y = \sigma_{yy}/\sigma_{zz}$ such that 
\[
	L_x=3|\alpha_x|,   \, \, \, L_y =3|\alpha_y| \, \, \,   \text { and } \,  \, \,  L_z=3 \, .
\]
In all computations the cell size is  constant $(5 \times 5 \times 5$ cells per unit volume) and is independent of the domain size.

\begin{figure}
  \centering \subfloat[][Positive stress ratios 
  $\sigma_{xx}/\sigma_{zz} > 0$ and $ \sigma_{yy}/\sigma_{zz} >  0$.]
  {\label{fig:matrixA}\includegraphics[scale=0.75]{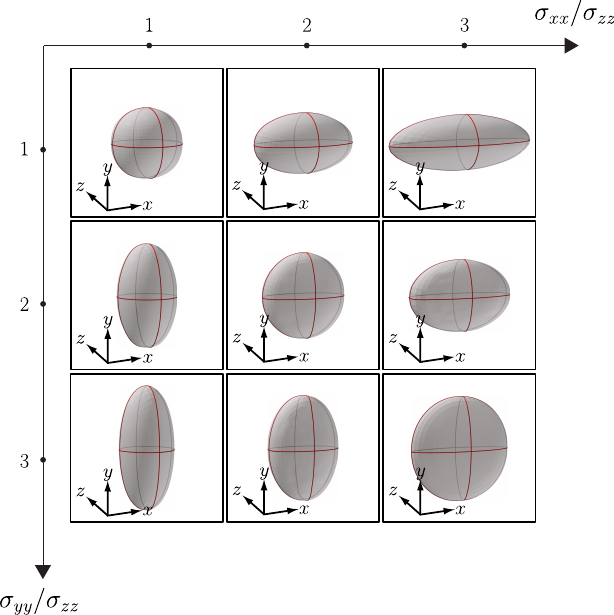}} \vspace{0.01\linewidth}
  \subfloat[][Negative stress ratios $\sigma_{xx}/\sigma_{zz}< 0$ and  $\sigma_{yy}/\sigma_{zz} < 0$.]
  {\label{fig:matrixB}\includegraphics[scale=0.75]{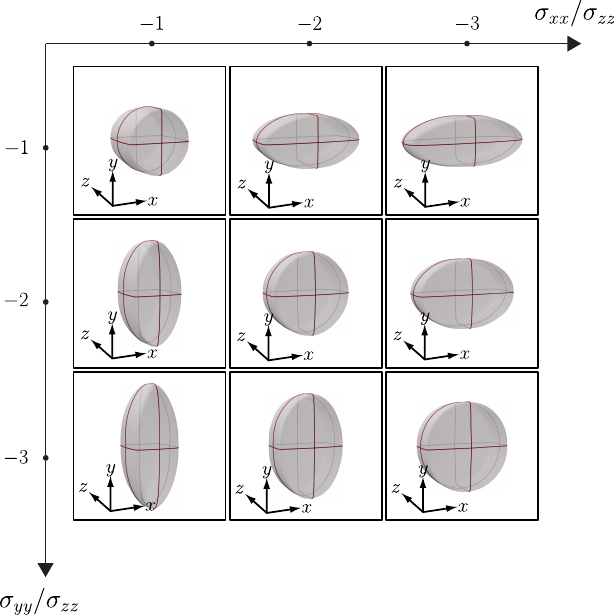}}
  \caption{Optimal hole shapes in a three-dimensional domain.  Computationally obtained hole shapes for different stress ratios.}
  \label{fig:optCavityShapesMatrix3D}
\end{figure}

As in the case of the preceding two-dimensional example in Section~\ref{sec:holesInPlate2D}, for negative stress ratios the compliance cost function~\eqref{eq:compliance} is augmented with an integral penalising surface area changes, cf.~\eqref{eq:compliancePenalty}.  The penalty modified cost function prohibits the formation of multiple holes, which cannot be obtained with shape optimisation.  The obtained hole shapes are shown in Figure~\ref{fig:optCavityShapesMatrix3D}. The applied positive stress ratios result in ellipsoidal hole shapes with semi-axis ratios proportional to the stress ratios (up to discretisation errors).  On the other hand, the applied negative stress ratios result  in hole shapes where the cross-section in the $x-y$ plane is an ellipse and in the  $x-z$ and $y-z$ planes  are smoothed quadrilaterals. 
%

\subsection{Shape and topology optimisation \label{sec:shapeTopOptEx}}
%
This section introduces several examples with combined shape and topology optimisation. The topology of the domain is altered by introducing new holes based on the topology derivative. Subsequently, the  shape of the holes  is optimised with the presented multiresolution shape optimisation technique. This approach is in spirit similar to the classical  \emph{bubble method}  by Eschenauer et al.~\cite{eschenauer1994bubble}.  In our present  implementation we do not consider the merging or removing of holes, hence our approach is  more restrictive than most newer topology optimisation techniques. For excellent recent reviews on various topology optimisation techniques see \cite{sigmund2013topology,van2013level,deaton2014survey}.  

Without going in to details, the topology derivative at a point gives the change of the cost function when a small hole is introduced at that point, see e.g. \cite{eschenauer2001topology, garreau2001topological, novotny2003topological}. To this end, in addition to the original domain $\Omega$ a modified domain $\Omega_r$ containing a hole of radius of $r$ is considered. The hole shape is  either a circle or sphere depending on the dimension of the domain.  The cost function $J(\Omega_r, \vec u_r) $ of the problem defined on $\Omega_r$ can be  expressed with a series expansion 
\begin{equation} \label{eq:topologySeries}
	J(\Omega_r, \vec u_r)  =  J(\Omega , \vec u ) + f (r) D_T J(\Omega , \vec u ) + \operatorname {O}( r^d) \, 
\end{equation}
where $f(r)$ is the size of the hole and  $D_T J(\Omega , \vec u )  $ is the topology derivative.  The hole size in 2D is  $f(r)= r^2 \pi$ and in 3D it is $f(r)= 4\pi r^3/3$. The topology derivative can be related to  the shape derivative by considering the expansion of an infinitesimally small hole~\cite{novotny2003topological}.  In our implementation we use the expressions given in \cite{novotny2006topological, novotny2007topological}   for the topology derivative of the compliance cost function. Depending on the dimension of the domain we obtain the topology derivative for a material with Poisson's ration $\nu$ with the following equations:
\begin{itemize}
\item two-dimensional elasticity, plane-stress
	\begin{equation}
		\begin{aligned}
		D_T  J (\Omega, \vec u) &= \frac{4}{1+\nu} \vec {\sigma} (\vec{u}) : \vec {\epsilon} (\vec{u}) \\ &- \frac{1-3\nu}{1-\nu^2} \trace \vec{\sigma} (\vec u) \trace 		\vec {\epsilon} (\vec u)	
		\end{aligned}
	\end{equation}
\item two-dimensional elasticity, plane-strain
	\begin{equation}
		\begin{aligned}
		D_T   J (\Omega, \vec u)  &=   4 (1-\nu) \vec {\sigma} (\vec{u}) : \vec {\epsilon} (\vec{u})  \\ &- \frac{(1-4\nu) (1-\nu) }{1-2 \nu }\trace \vec{\sigma} (\vec u) \trace \vec {\epsilon} (\vec u)	
		\end{aligned}
	\end{equation}
\item three-dimensional elasticity
	\begin{equation}
		\begin{aligned}
		D_T ( J (\Omega, \vec u)  = & \frac{3}{2}  \frac{1-\nu}{7-5\nu}  \bigg [  10  \vec \sigma(\vec u) : \vec \epsilon (\vec u )   \\ & - \frac{1-5\nu}{1-2\nu}  \trace \vec 	\sigma(\vec u)   \trace \vec \epsilon (\vec u)	 \bigg  ]
		\end{aligned}
	\end{equation}
\end{itemize}

We use the isocontours of the topology derivative  $D_T (\vec x) $  to determine the location and shape of the hole to be introduced. To this end, we introduce a control mesh which approximates the boundary of the  region to be removed from the domain. Although this step is presently performed manually, it is feasible to automate it. This is particularly straightforward in case of two-dimensional domains with a polygon as the boundary of the hole. For three-dimensional problems the  isocontour of the topology derivative can be first extracted with  a marching-cube algorithm and subsequently remeshed in order obtain a  uniform control mesh  with mostly regular vertices, see e.g.~\cite{botsch2010polygon}. 

After a hole is generated, it is subjected to shape optimisation. During the shape optimisation the size of the hole is allowed to increase by a prescribed amount. For instance, the area constraint~\eqref{eq:areaConst} is modified as follows
\begin{equation}   \label{eq:increaseArea} 
	A(\vec x^{\ell_c}) \ge \rho_A A_0 \, , 
\end{equation}
where $\rho_A$ is a user prescribed scalar and $A_0$ is the area of the initial hole. In this context the boundary of each hole is treated as a separate multiresolution curve or surface.

In passing we note that instead of the topology derivative it would also be possible to determine the location and shape of the introduced holes with the density based SIMP technique widely used in topology optimisation~\cite{bendsoe1989optimal, bendsoe2003topology}.  

\subsubsection{Cantilever \label{sec:cantilverExample}}
\begin{figure}[!th]
  \centering 
  \subfloat[][Problem description. The top boundary is  first optimised with constraints on the end nodes (solid squares) to prevent  horizontal movement. ] {\label{fig:trussTopoA}
   \includegraphics[scale=1]{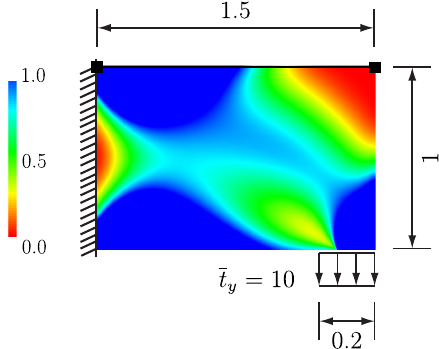}}\\
  \vspace{0.01\textheight} 
  \subfloat[][Two holes are introduced at locations of minimum topology derivative.] {\label{fig:trussTopoB}
   	\includegraphics[scale=1]{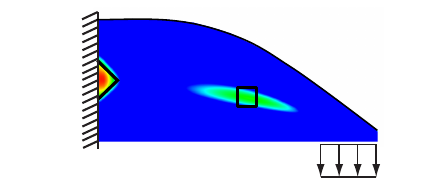}}  \\ 
	\vspace{0.01\textheight} 
  \subfloat[][Final optimised geometry.]{\label{fig:trussTopoC}
    \includegraphics[scale=1]{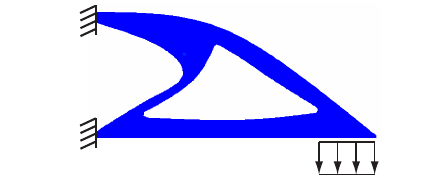}}
  \caption{Shape and topology  optimisation of a cantilever truss.  The isocontours indicate the topology derivative.}
  \label{fig:trussTopo}
\end{figure}

%
The cantilever shown in Figure~\ref{fig:trussTopoA} is a widely studied benchmark example in topology and shape optimisation. Most relevant to our study are the results  obtained by  Eschenauer et al.~\cite{eschenauer1994bubble} using the bubble method, which is as previously mentioned is similar to our approach. In our computations the rectangular domain is loaded with a  distributed force $\overline{t}_y = 10 $ close to the lower right corner and  the Young's modulus and Poisson's ratio are chosen with $E=100$ and $\nu=0.4$, respectively.  The immersed finite element grid contains $100 \times 100$ cells of uniform size. 

The final topology and shape optimised geometry is shown in Figure~\ref{fig:trussTopoC}. This geometry is similar to the ones presented in~\cite{eschenauer1994bubble, wang2003level} and has been obtained in three steps. 
\begin{enumerate}
	\item[1.]  The first step is a shape optimisation step.  During optimisation only the top boundary of the plate is allowed to move based on the computed shape gradients.  In the optimisation level $\ell_o=0$ the top boundary is represented with a control polygon using $2$ elements and $3$ vertices.  In the twice subdivided computational level $\ell_c=2$  the control polygon  contains $16$ elements. During subdivision the two end nodes are tagged as corner so that they do not move horizontally.  In addition, during iterative shape optimisation the right corner node is allowed to move only vertically.  Moreover, we apply area  constraints of the form~\eqref{eq:increaseArea} so that the domain size is reduced. The geometry obtained with shape optimisation is shown in Figure~\ref{fig:trussTopoB}.
	\item[2.] The second step is a topology optimisation step. 	After computing the topology derivative we manually introduce two holes at locations with minimum topology derivative, see Figure ~\ref{fig:trussTopoB}.  The first hole is triangular shaped and splits the clamped boundary into an upper and lower part. The second hole is square shaped and is located inside the domain. The holes have to be large enough so that they can be represented on the immersed finite element grid. According to~\cite{Ruberg:2014aa} the minimum hole size has to be larger than 
\begin{equation} 
	2 \sqrt{d} (\alpha +1) h \, ,
\end{equation}	
where $d=2$ is the dimensionality of the domain, $\alpha$ is the polynomial degree of the b-spline basis functions and $h$ is the cell size.
	\item[3.]	 The last step is again a shape optimisation step. The control polygons belonging to the previously introduced two holes with three and four nodes, respectively, are subdivided twice to obtain the computational control polygon at level  $\ell_c=2$.  Subsequently the geometry of the two holes is  iteratively optimised using shape optimisation. The optimisation is terminated before the two holes start to merge.
\end{enumerate}

%
\subsubsection{Simply supported plates with different aspect ratios}
\begin{figure}[ht]
  \centering 
  {\includegraphics[scale=1]{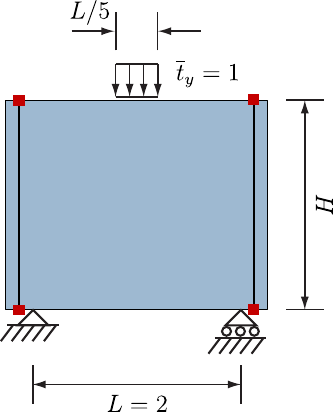}}
  \caption[Arch to truss transition, problem description]{Simply supported plate with aspect ratio $H/L$. Each of the two supports have a width of $0.05L$. The two bold vertical lines indicate the free boundaries which are allowed to move during shape optimisation. The two ends of each of the vertical lines (red squares) are constraint to move only horizontally.}
  \label{fig:archToTrussDescription}
\end{figure}

%
We consider three simply supported plates with different aspect ratios $H/L$, see Figure \ref{fig:archToTrussDescription}. It is known from structural analysis that for slender plates a truss-like system and for stocky plates an arch-like system is more efficient. The plate is loaded with a symmetrically placed distributed vertical force $\overline{t}_y =1$ of length $L/5$ on its top edge. The Young's modulus and Poisson's ratio are chosen with $E=100$ and $\nu=0.4$, respectively. We consider three different aspect ratios $H/L \in \{0.25,0.5,1\}$ with the corresponding immersed finite element grids  containing  $80\times200$, $200\times100$ and $150\times150$ cells, respectively. In the computations the height is fixed to $H=1$ and only the length $L$ is varied.

\begin{figure*}
  \centering 
  \includegraphics[width=0.9\textwidth]{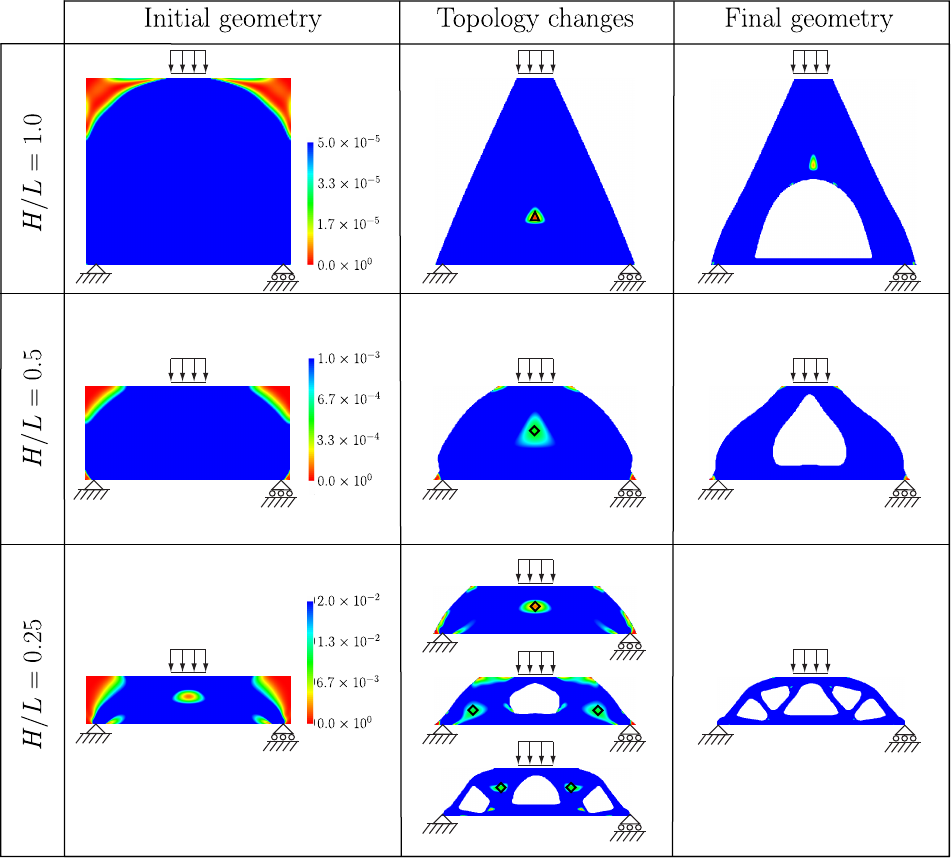}
  \caption{Simply supported plates with different aspect ratios.  In all snapshots the isocontours indicate the topology derivative. In the middle column the small triangular and square shaped holes introduced during the topology optimisation step are shown.}
  \label{fig:archToTrussFinalTopo}
\end{figure*}

As in the cantilever example, cf. Section~\ref{sec:cantilverExample},  the final geometry is obtained in several steps, namely an initial shape optimisation step is followed by several topology optimisation and shape optimisation steps. In the initial shape optimisation step the two vertical boundaries of the plate are optimised while the domain area is allowed to reduce. At the coarse optimisation level $\ell_o=0$ each edge is represented with two elements, which are three times subdivided to obtain the computational control polygon at level $\ell_c=3$. In the subsequent topology optimisation step we semi-manually introduce triangle and square shaped holes, see Figure~\ref{fig:archToTrussFinalTopo} middle column. For each of the polygons the optimisation and computation levels are chosen with $\ell_o=0$ and $\ell_c=3$, respectively.  The topology optimisation is followed by the shape optimisation of all the domain boundaries. The movement of the vertices under the distributed force and supports are constrained to remain fixed.  In case of  the plate with small  $H/L=0.25$ several more topology and shape optimisation steps are performed, see Figure~\ref{fig:archToTrussFinalTopo} middle column. The obtained geometries are shown in Figure~\ref{fig:archToTrussFinalTopo} right column. The optimisation is always terminated before any holes start to merge. As expected, we obtain for a $H/L=1.0$ an arch-like structure and for  $H/L=0.25$ a truss-like structure. 

%

\subsubsection{Three-dimensional stool}
\begin{figure}
  \centering 
  {\includegraphics[scale=1]{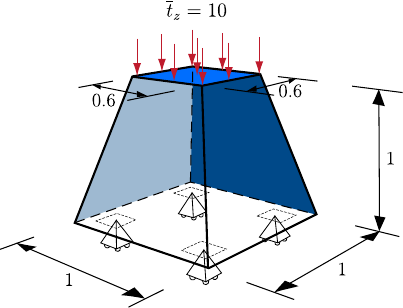}}
  \caption{Three-dimensional stool. Problem description.  Roller supports are applied to all finite element nodes inside the regions of size $0.2 \times 0.2$ marked by dashed squares.   \label{fig:tableDomain}}
\end{figure}

%
In this last example we present the combined topology and shape optimisation of a three-dimensional solid, see Figure~\ref{fig:tableDomain}. The initial domain is a truncated pyramid and is at its top loaded with a uniform distributed load  $\overline{t}_z=10$.  At its bottom it is supported by four distributed roller supports each of size $0.2 \times 0.2$.  The Young's modulus and Poisson's ratio are chosen with  $E=100$ and $\nu=0.4$, respectively.

In the optimisation study only one quarter of the domain is considered and appropriate bounds and geometry tags are applied at the two planes of symmetry. The corresponding immersed finite element grid is of size  $0.7\times0.7\times1$ and consists of $30\times30\times30$ cells. 

\begin{figure*}
  \centering 
  	\subfloat[][Initial geometry.]{
	\includegraphics[scale=1]{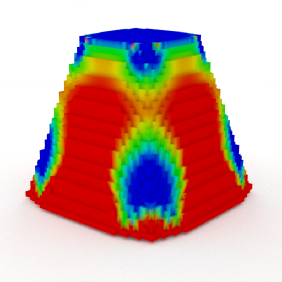}}
	  \hspace{0.1\linewidth} 
	\subfloat[][Optimised geometry. \label{fig:firstToptOptResult}]{
   	\includegraphics[scale=1]{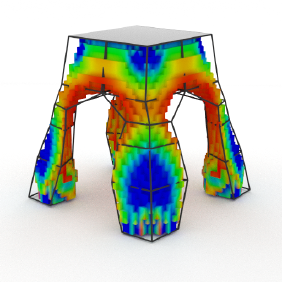}}\\
     \caption{Three-dimensional stool. First topology optimisation step. The isocontours indicate the topology derivative. The wireframe represents the coarse control mesh.}
	  \label{fig:tableOptimisationFirstTopOpt}
\end{figure*}

\begin{figure*}
  \centering
	\subfloat[][Initial geometry.]{
   	\includegraphics[scale=1]{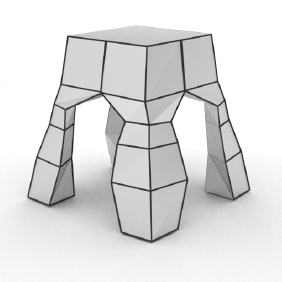}}
	\hspace{0.1\linewidth} 
	\subfloat[][Optimised geometry.]{
    	\includegraphics[scale=1]{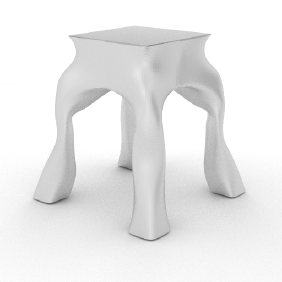}}
  	\caption{Three-dimensional stool. First shape optimisation step.  The wireframe represents the coarse resolution mesh used for optimisation.}
  \label{fig:tableOptimisationFirstShapeOpt}
\end{figure*}

\begin{figure*}
	\centering 
	\subfloat[][Initial geometry.]{
	\includegraphics[scale=1]{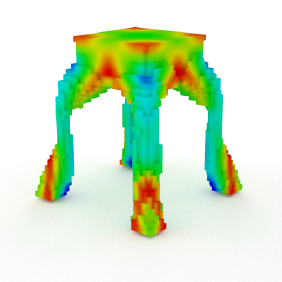}}
  	\hspace{0.1\linewidth} 
	\subfloat[][Optimised geometry. \label{fig:secondToptOptResult}]{
	\includegraphics[scale=1]{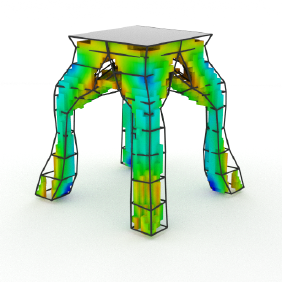}}
    \caption{Three-dimensional stool. Second topology optimisation step. The isocontours indicate the topology derivative. The wireframe represents the coarse control mesh.}
  \label{fig:tableOptimisationSecTopOpt}
\end{figure*}

\begin{figure*}
	\centering
	\subfloat[][Initial geometry.]{
    	\includegraphics[scale=1]{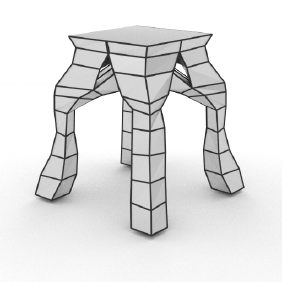}}
  	\hspace{0.1\linewidth} 
	\subfloat[][Optimised geometry.]{
 	\includegraphics[scale=1]{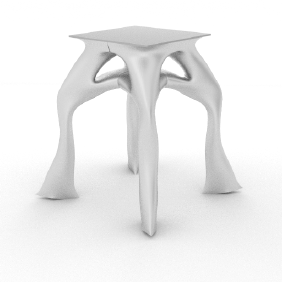}}
 	\caption{Three-dimensional stool. Second shape optimisation step.  The wireframe represents the coarse resolution mesh used for optimisation.}
   \label{fig:tableOptimisationSecShapeOpt}
\end{figure*}

The sequence of  the performed topology and shape optimisation steps are shown in Figures~\ref{fig:tableOptimisationFirstTopOpt}, \ref{fig:tableOptimisationFirstShapeOpt}, \ref{fig:tableOptimisationSecTopOpt} and \ref{fig:tableOptimisationSecShapeOpt}. In total two topology and two shape optimisation steps are performed. In each topology optimisation step we remove  in one go a relatively large amount of material by deleting computational cells with topology derivative below a threshold. In the first topology optimisation step illustrated in Figure~\ref{fig:tableOptimisationFirstTopOpt} all cells with topology derivative $D_T  J (\Omega, \vec u) \leq 0.025$ are removed. Subsequently, we semi-manually generate the coarse resolution subdivision control mesh depicted in Figure~\ref{fig:firstToptOptResult} for representing the new topology.  In the following shape optimisation step, see Figure~\ref{fig:tableOptimisationFirstShapeOpt},   the generated control mesh serves as the  optimisation level $\ell_o=0$ and the computation level is chosen with $\ell_c=2$. During the shape optimisation the volume of the domain is constraint to remain constant. In the second topology optimisation step shown in Figure~\ref{fig:tableOptimisationSecTopOpt}  all cells with topology derivative $D_T J (\Omega, \vec u) \leq 0.04$ are removed. This is followed by a semi-manual control mesh generation, see Figure~\ref{fig:secondToptOptResult}, and the final shape optimisation step shown Figure~\ref{fig:tableOptimisationSecShapeOpt}.

%

\section{Summary and conclusions}

A multiresolution optimisation technique  based on subdivision surfaces was introduced. The domain boundaries are described with subdivision surfaces and  the domain boundary value problem is discretised with an immersed finite element technique. The wavelet-like multiresolution decomposition of the domain boundary  yields a low resolution control mesh and a sequence of detail vectors. The control mesh at any specific refinement level can be reconstructed  on-the-fly with the introduced synthesis operator. Editing  the coarse levels leads to large-scale geometry changes while editing fine levels leads to small-scale geometry changes. In addition, the multiresolution editing semantics allows the decoupling of the choice of the editing level from the  size of the geometric features present in the geometry. In the proposed approach we start optimising the coarsest control mesh and successively increase the optimisation level each time a minimum is reached.  The domain geometry is always described with a fine control mesh on the immersed finite element grid,  independently from the control mesh level used for optimisation. Hence, any fine scale geometric details, like fillets or surface undulations,  are faithfully represented on the discretisation grid. As our  examples demonstrate multiresolution shape optimisation enables us to find better optima and is exceedingly robust, partly due to the use of the immersed finite element technique.

The multiresolution editing semantics appears to be particularly appealing for isogeometric analysis because it enables the full decoupling of the geometry and the analysis representations of the same geometry. It allows to seamlessly map variables and fields between the two representations irrespective of their resolutions.  Beyond optimisation this decoupling can be useful in a number of applications, such as for computing fast approximate solutions or  multigrid and multilevel preconditioners~\cite{Borzi2013}. The presented multiresolution techniques can also be extended to non-uniform rational b-splines (NURBS), which are the more commonly used basis functions in industrial software.  It is straightforward to include rational b-splines through the use of homogenised coordinates, see, e.g.,~\cite{Cirak:2011aa}. In order to consider non-uniform b-splines it is instructive to consult previous works on non-uniform subdivision~\cite{Sederberg:2003aa, Cashman:2009aa}. Although we focused in the present paper on uniform subdivision refinement and coarsening, it is conceivable (and desirable) to develop adaptive multiresolution algorithms in the spirit of hierarchical b-splines~\cite{Grinspun:2002aa, Bornemann:2013aa,wei2015truncated}.  The utility of adaptive geometry representations in shape optimisation has already been demonstrated, for instance, with b-spline surfaces refined by knot insertion~\cite{han2014adaptive} and reparameterisation of Bezier surfaces~\cite{desideri2007nested}. Furthermore, in the present paper in 3D holes were introduced by manually fitting control meshes to the isocontours of the topology derivative. This process can be automated using techniques for extracting subdivision control meshes from isocontours~\cite{bommes2013quad}.  The related issue of merging of holes can be achieved with approximate Boolean operations for subdivision surfaces~\cite{Biermann:2001aa}. In closing, it is worth emphasising that basic subdivision techniques recently became available in a number of engineering design software, including Autodesk Fusion 360, PTC Creo and CATIA, which most likely will increase their use  in future engineering practice.

\section*{Acknowledgement}
The partial support of the EPSRC through grant \# EP/G008531/1  and EC through Marie Curie Actions (IAPP)  program CASOPT project are gratefully
acknowledged.

\bibliographystyle{elsarticle-num-names}
\bibliography{./bloxOptimisation}

\end{document}